\documentclass{article}

\usepackage{amsmath}
\usepackage{amsfonts}
\usepackage{amssymb}
\usepackage{xspace}
\usepackage{graphicx}
\makeatletter \@addtoreset{equation}{section} \makeatother

\numberwithin{equation}{section}
\makeatletter \@addtoreset{equation}{section} \makeatother

\newcommand{\adl}{\vspace{1\baselineskip}}

\newenvironment{proof}{\par\noindent{\sc Proof:}
}{\hfill\llap{$\Box$}\vspace{1\baselineskip}\par\noindent}
\newenvironment{proofof}{\par\noindent{\sc Proof}
}{\hfill\llap{$\Box$}\vspace{1\baselineskip}\par\noindent}
\newtheorem{theorem}{Theorem}[section]
\newtheorem{proposition}[theorem]{Proposition}
\newtheorem{lemma}[theorem]{Lemma}
\newtheorem{corollary}[theorem]{Corollary}
\newtheorem{remark}[theorem]{Remark}
\newtheorem{definition}[theorem]{Definition}
\newtheorem{example}[theorem]{Example}

\newcommand{\beq}{\begin{equation}}
\newcommand{\eeq}{\end{equation}}

\newcommand{\ba}{\begin{array}}
\newcommand{\ea}{\end{array}}
\newcommand{\bt}{\begin{theorem}}
\newcommand{\et}{\end{theorem}}
\newcommand{\bp}{\begin{proposition}}
\newcommand{\ep}{\end{proposition}}
\newcommand{\bl}{\begin{lemma}}
\newcommand{\el}{\end{lemma}}
\newcommand{\bc}{\begin{corollary}}
\newcommand{\ec}{\end{corollary}}
\newcommand{\bi}{\begin{itemize}}
\newcommand{\ei}{\end{itemize}}
\newcommand{\ben}{\begin{enumerate}}
\newcommand{\een}{\end{enumerate}}
\newcommand{\bpf}{\begin{proof}}
\newcommand{\epf}{\end{proof}}
\newcommand{\bpff}{\begin{proofof}}
\newcommand{\epff}{\end{proofof}}

\newcommand{\bdf}{\begin{definition}\rm}
\newcommand{\edf}{\end{definition}}
\newcommand{\br}{\begin{remark}\rm}
\newcommand{\er}{\end{remark}}
\newcommand{\bex}{\begin{example}\rm}
\newcommand{\eex}{\end{example}}
\def\pri{\hbox to 10pt{\hfil\hbox to 0.4pt{\vrule height5pt width0.4pt
                 depth0pt}\vrule width5pt height0.4pt depth0pt\hfil}}
\newcommand{\TC}{{\rm TC}}
\newcommand{\CT}{{\rm CT}}

\newcommand{\TAT}{{\rm TAT}}

\newcommand{\gk}{{\bf{k}}}

\newcommand{\gK}{{\bf K}}

\newcommand{\gT}{{\bf T}}

\newcommand{\calL}{{\cal L}}
\newcommand{\calK}{{\cal K}}
\newcommand{\calT}{{\cal T}}

\newcommand{\gR}{{\mathbb R}}

\newcommand{\gN}{{\mathbf N}}
\newcommand{\gU}{{\mathbf U}}

\newcommand{\gZ}{{\mathbb Z}}

\newcommand{\Nat}{{\mathbb N}}
\newcommand{\Rn}{{\mathbb R}^{n}}
\newcommand{\RN}{{\mathbb R}^{N}}

\newcommand{\D}{{\cal D}}

\newcommand{\Ha}{{\cal H}}

\newcommand{\BV}{\mathop{\rm BV}\nolimits}

\def\sgn{\mathop{\rm sgn}\nolimits}

\newcommand{\If}{\ \mbox{\rm if}\ }

\newcommand{\wc}{\rightharpoonup}

\newcommand{\Sph}{{\mathbb S}}
\newcommand{\RP}{{\mathbb{RP}}^2}
\newcommand{\wRP}{{\rm RP^{2}}}
\newcommand{\SP}{{\Sph^2}}

\def\mesh{\mathop{\rm mesh}\nolimits}
\newcommand{\ttt}{{\mathfrak t}}
\newcommand{\nnn}{{\mathfrak n}}
\newcommand{\bbb}{{\mathfrak b}}
\newcommand{\kk}{{\mathfrak K}}
\newcommand{\gtau}{{\boldsymbol{\tau}}}

\newcommand{\Var}{\mathop{\rm Var}\nolimits}

\let\a=\alpha
\let\be=\beta

\let\d=\delta
\let\e=\varepsilon

\let\g=\gamma

\let\l=\lambda
\let\m=\mu
\let\n=\nu

\let\p=\pi

\let\r=\rho
\let\s=\sigma
\let\t=\theta
\let\tt=\tau
\let\x=\xi

\let\GG=\Gamma

\let\wih=\widehat

\let\wid=\widetilde

\let\sb=\subset

\let\Lra=\Longrightarrow

\let\fa=\forall
\let\tim=\times

\let\ul=\underline
\let\ds=\displaystyle

\let\lan=\langle
\let\ran=\rangle

\let\i=\infty

\title{\Large \bf The weak Frenet frame of non-smooth curves \\ with finite total curvature and absolute torsion}
\author{\it Domenico Mucci and Alberto Saracco
\footnote{%
{\sc Dipartimento di Scienze Matematiche,
Fisiche ed Informatiche, Universit\`{a} di Parma,
Parco Area delle Scienze 53/A, I-43124 Parma, Italy.
E-mail: domenico.mucci@unipr.it, alberto.saracco@unipr.it}
}
}
\begin{document}
\date{}
%\date{Version of \today}
\topskip=1.5truecm \maketitle \topskip=1.5truecm \maketitle
       %%%%%%%%% Abstract %%%%%%%%%%%%%%%%
%
\begin{abstract} We deal with a notion of weak binormal and weak principal normal for non-smooth curves of the Euclidean space with finite total curvature and total absolute torsion.
By means of piecewise linear methods, we first introduce the analogous notion for polygonal curves, where the polarity property is exploited, and then make use of a density argument.
Both our weak binormal and normal are rectifiable curves which naturally live in the projective plane.
In particular, the length of the weak binormal agrees with the total absolute torsion of the given curve.
Moreover, the weak normal is the vector product of suitable parameterizations of the tangent indicatrix and of the weak binormal.
In the case of smooth curves, % with positive curvature,
the weak binormal and normal yield (up to a lifting) the classical notions of binormal and normal.
Finally, the torsion force is introduced: similarly as for the curvature force, it is a finite measure obtained by performing the tangential variation of the length of the tangent indicatrix in the Gauss sphere.
\end{abstract}

\adl\par


\section{Introduction}
In classical differential geometry, it sometimes happens that the geometry of a proof can become obscured by analysis.
This statement by M.~A.~Penna \cite{Pe}, which may be referred e.g. to the classical proof of the Gauss-Bonnet theorem, suggests to apply piecewise linear methods in order to make the geometry of a proof completely transparent.
\par For this purpose, by using the geometric description of the torsion of a smooth curve, Penna \cite{Pe} gave in 1980 a suitable definition of {\em torsion}
for a polygonal curve of the Euclidean space $\gR^3$, and used piecewise linear methods and homotopy arguments to produce an illustrative proof of the well-known property that the total torsion of any closed unit speed regular curve of the unit sphere $\SP$ is equal to zero.
\par Differently to the smooth case, the polygonal torsion is a function of the segments.
His definition, in fact, relies on the notion of {\em binormal vector} at the interior vertices.
Since the angle between consecutive discrete binormals describes the movements of the ``discrete osculating planes" of the polygonal, binormal vectors naturally live in the {\em projective plane $\RP$}, see Sec.~\ref{Sec:poly}.
\par We recall here that J.~W. Milnor~\cite{Mi} defined the {\em tangent indicatrix}, or tantrix, of a polygonal $P$ as the geodesic polygonal $\ttt_P$ of the Gauss sphere $\SP$
obtained by connecting with oriented geodesic arcs the consecutive points given by the direction of the oriented segments.
Therefore, the {\em total curvature} $\TC(P)$, i.e., the sum of the turning angles of the polygonal, agrees with the length $\calL_{\SP}(\ttt_P)$ of the tantrix, and the
total absolute torsion $\TAT(P)$ agrees with the sum of the shortest angles
between the geodesic arcs in $\SP$ meeting at the edges of $\ttt_P$.
%
%, i.e., with the total curvature of the tantrix in $\SP$.
%Of course, the two above definitions of total absolute torsion are equivalent, compare Remark~\ref{Requal}.
%
\par From another viewpoint, W.~Fenchel \cite{Fe} in the 1950's exploited the {\em spherical polarity} of the tangent and binormal indicatrix in order to analyze the differential geometric properties of smooth curves in $\gR^3$. In his survey, Fenchel proposed a general method that gathers several results on curves in a unified scheme.
We point out that Fenchel deals with $C^4$ rectifiable curves (parameterized by arc-length) such that at each point it is well-defined the osculating plane,
that is, a plane containing the linearly independent vectors $\ttt:=\dot c$ and $\ddot c$, such that its suitably oriented normal unit vector $\bbb$, the binormal vector, is of class $C^2$, and the two vectors $\dot \ttt$ and $\dot\bbb$ never vanish simultaneously. He then defines the principal normal by the vector product
\beq\label{int-normal}
\nnn:=\bbb\tim\ttt\,.  \eeq
Since the derivatives of $\ttt$ and $\bbb$ are perpendicular to both $\ttt$ and $\bbb$, the curvature $\gk$ and torsion $\gtau$ are well-defined through the formulas:
$$ \dot\ttt=\gk\,\nnn\,,\qquad \dot\bbb=-\gtau\,\nnn\,.$$
As a consequence, one has
$$ \dot\nnn=-\gk\,\ttt+\gtau\,\bbb $$
and hence the Frenet-Serret formulas hold true, but Fenchel allows both curvature and torsion to be zero or negative.
Related arguments have been treated in \cite{Ba,EI,EIS,Mc,To}.
\adl\par\noindent{\large\sc Content of the paper.} We deal with curves in the Euclidean space $\gR^3$ with finite total curvature and total absolute torsion.
We address to J.~M.~Sullivan \cite{Su_curv} for the analysis of curves with finite total curvature, and also to our paper \cite{MuSa} for the $\BV$-properties of the unit normal of planar curves.
\par By melting together the approaches by Penna and Fenchel previously described, in this paper we firstly define the {\em binormal indicatrix} $\bbb_P$ of a polygonal $P$ in $\gR^3$ as
the arc-length parameterization $\bbb_P$ of the polar in $\RP$ of the tangent indicatrix $\ttt_P$, see Definition~\ref{Dbin} and Figure~\ref{fig:SferaTB}.
Therefore, the {\em total absolute torsion} $\TAT(P)$ of $P$ is equal to the length of the curve $\bbb_P$.
We remark that a similar definition has been introduced by T.~F.~Banchoff in his paper \cite{Ba} on space polygons.
\par However, differently from what happens for the length and the {\em total curvature}, the monotonicity formula fails to hold. More precisely, if $P'$ is a polygonal inscribed in $P$, by the triangular inequality we have $\calL(P')\leq \calL(P)$ and $\TC(P')\leq \TC(P)$, compare e.g. \cite[Cor.~2.2]{Su_curv}, but it may happen that $\TAT(P')>\TAT(P)$,
see Example~\ref{Emon}.
\par For that reason, the {\em total absolute torsion} $\TAT(c)$ of a curve $c$ in $\gR^3$ is defined by following the approach due to Alexandrov-Reshetnyak \cite{AR},
that involves the notion of {\em modulus} $\m_c(P)$ of a polygonal $P$ inscribed in $c$, see \eqref{DTAT}.
\par As a consequence, by means of a density argument, a good notion of {\em weak binormal indicatrix} $\bbb_c$ for a non-smooth curve with {\em finite total curvature} and
{\em absolute torsion} is obtained in our first main result, Theorem~\ref{Tbin}.
In fact, we infer that for any sequence $\{P_h\}$ of inscribed polygonals with $\m_c(P_h)\to 0$, one has $\TAT(P_h)\to\TAT(c)$, see Proposition~\ref{PTAT},
and hence that the weak binormal $\bbb_c$ only depends on the curve $c$.
%
%\par For this purpose, we recall that as for the length $\calL(c)$, the total curvature $\TC(c)$ and the total absolute torsion $\TAT(c)$ of a curve $c$ in $\gR^3$ are defined as the supremum of the corresponding value obtained among the inscribed polygonals. Moreover, they can be recovered in terms of any sequence of inscribed polygonals with infinitesimal meshes,
%compare Sullivan \cite{Su_curv} and Proposition~\ref{PTAT}.
%
\par For smooth curves, the total absolute torsion, which agrees with the length in the Gauss sphere of the smooth binormal curve $\bbb$,
actually agrees with the total geodesic curvature of the smooth tantrix $\ttt$ in $\SP$.
\par In fact, on account of the density result from \cite[Prop.~4]{Pe}, by Proposition~\ref{PTAT} one readily obtains that
\beq\label{int-TAT} \TAT(c)=\int_c |\gtau|\,ds \eeq
where $\gtau(s)$ is the torsion of the smooth curve $c$.
This property is checked in Example~\ref{Ehelix}, referring to a helicoidal curve, where we exploit piecewise linear methods in the computation.
\smallskip\par In Theorem~\ref{Tbin}, we show the existence of a curve $\bbb_c$ of $\RP$, parameterized by arc-length, whose length is equal to the total absolute torsion:
\beq\label{int-Lbc} \calL_{\RP}(\bbb_c)=\TAT(c)\,. \eeq
\par The hypothesis $\TC(c)<\i$ in Theorem~\ref{Tbin} may sound a bit unnatural, and actually a technical point, since it allows us to prove that $\bbb_c$ has constant velocity one, so that \eqref{int-Lbc} holds true.
\par To this purpose, we recall that the definition of {\em complete torsion} $\CT(P)$ of polygonals $P$ given by Alexandrov-Reshetnyak \cite{AR},
who essentially take the distance in $\SP$ between consecutive discrete binormals, implies that planar polygonals may have positive torsion at ``inflections points".
Defining the {\em complete torsion} $\CT(c)$ of curves $c$ in $\gR^3$ as the supremum of the complete torsion of the inscribed polygonals,
they obtain in \cite[p.~244]{AR} that any curve with finite complete torsion and with no points of return must have finite total curvature.
\par With our definition of torsion, the above implication clearly fails to hold, see Remark~\ref{RFrenet}.
On the other hand, equality \eqref{int-TAT} is violated if one considers the complete torsion from \cite{AR}, since
for a smooth planar curve with inflection points, one has $\CT(c)>0$.
\par We finally notice that a curve with finite total curvature and total absolute torsion may have infinite complete torsion in the sense of \cite{AR}:
just take a smooth planar curve with a countable set of inflection points.
%\par We finally observe that a curve $c$ has finite complete torsion in the sense of \cite{AR}, if and only if $\TC(c)+\TAT(c)<\i$.
%In fact, for polygonals $P$ one has $$ \max\{\TC(P),\,\TAT(P) \}\leq \CT(P)\leq \TC(P)+\TAT(P)\,. $$
%
%\smallskip
\par In Theorem~\ref{Tbinsmooth}, we show that for smooth curves whose torsion $\gtau$ (almost) never vanishes, the {\em weak binormal} $\bbb_c$ obtained in Theorem~\ref{Tbin}, when lifted to $\SP$,
agrees with the arc-length parameterization of the smooth binormal $\bbb$.
\par Similar features concerning the tantrix are collected in Propositions~\ref{Ptan} and
\ref{Ptansmooth}. Our curve $\ttt_c$ satisfies $\calL_{\SP}(\ttt_c)=\TC(c)$ and hence it is strictly related with the
{\em complete tangent indicatrix} in the sense of Alexandrov-Reshetnyak \cite{AR}.
\smallskip\par Now, when looking for a possible weak notion of principal normal, a drawback appears. In fact, in Penna's approach \cite{Pe}, the curvature of an open polygonal $P$ is a non-negative measure $\m_P$ concentrated at the interior vertices, whereas the torsion is a signed measure $\n_P$ concentrated at the interior segments, see Remark~\ref{Rmeas}.
Since these two measures are mutually singular, in principle there is no way to extend Fenchel's formula \eqref{int-normal} in order to define the principal normal.
\par To overcome this problem, in Sec.~\ref{Sec:nor} we proceed as follows.
Firstly, we choose two suitable curves
$ \wid\ttt_P,\wid\bbb_P:[0,C+T]\to\RP$, where $C=\TC(P)$ and $T=\TAT(P)$,
that on one side inherit the properties of the tangent and binormal indicatrix $\ttt_P$ and $\bbb_P$, respectively, and on the other side take account of the order in which curvature and torsion are defined along $P$.
More precisely, one of the two curves is constant when the other one parameterizes a geodesic arc, whose length is equal to the curvature or to the (absolute value of the) torsion at one vertex or segment of $P$, respectively.
As in Fenchel's approach, by exploiting the polarity of the curves $\wid\ttt_P$ and $\wid\bbb_P$, the weak normal of the polygonal is well-defined by the inner product
$$\nnn_P(s):=\wid \bbb_P(s)\tim \wid\ttt_P(s)\in\RP\,,\qquad  s\in[0,T+C] $$
compare Remark~\ref{Rprod} and Figure~\ref{fig:SferaN}. Notice that by our Definition~\ref{Dnor} we have:
$$\calL_{\RP}(\nnn_P)=\TC(P)+\TAT(P)\,. $$
\par As a consequence, in our second main result, Theorem~\ref{Tnor}, using again an approximation procedure, the {\em weak principal normal} of a curve $c$ with finite total curvature and finite complete torsion is well-defined as a rectifiable curve $\nnn_c$ in $\RP$.
We recall that condition $\CT(c)<\i$ is stronger than the more natural assumption $\TAT(c)<\i$.
It turns out that the product formula \eqref{int-normal} continues to hold in a suitable sense, and we also have:
$$ \calL_{\RP}(\nnn_c)=\TC(c)+\TAT(c)\,.$$
\par In particular, for smooth curves whose curvature (almost) never vanishes, the principal normal $\nnn$ agrees with a lifting of a suitable parameterization of the weak normal $\nnn_c$.
More precisely, in Proposition~\ref{Pnorsmooth} we obtain that
$$[\nnn(s(t))]=\nnn_c(t)\in\RP \qquad \fa\, t\in[0,\TC(c)+\TAT(c)] $$
where $s(t)$ is the inverse of the increasing and bijective function
$$ t(s):=\int_0^s(\gk(\l)+|\gtau(\l)|)\,d\l\,,\qquad s\in[0,\calL(c)] \,. $$
\par In Sec.~\ref{Sec:smooth}, we make use of an analytical approach in order to define the binormal and principal normal
of smooth regular curves with inflection points. In fact, if $|\dot c(s_0)|=1$ but $\ddot c(s_0)=0_{\gR^3}$, the first non-zero higher order derivative $c^{(n)}(s_0)$ of $c$ at $s_0$ satisfies $\dot c(s_0)\perp c^{(n)}(s_0)$ and hence it plays a role in the definition of the binormal. Therefore, following Fenchel \cite{Fe} in order to define the normal as in \eqref{int-normal},
in Proposition~\ref{PFS} we get:
$$ \ttt(s_0)=\dot c(s_0)\,,\quad \bbb(s_0)=%{\frac{\dot c(s_0)\tim c^{(n)}(s_0)}{ |\dot c(s_0)\tim c^{(n)}(s_0)|}}=
{\frac{\dot c(s_0)\tim c^{(n)}(s_0)}{ \Vert\dot c^{(n)}(s_0)\Vert}}\,,\quad
\nnn(s_0)=\frac{c^{(n)}(s_0) }{ \Vert c^{(n)}(s_0)\Vert}\,. $$
%
%
%$$ \bbb(s_0)={\frac{\dot c(s_0)\tim c^{(n)}(s_0)}{ |\dot c(s_0)\tim c^{(n)}(s_0)|}}={\frac{\dot c(s_0)\tim c^{(n)}(s_0)}{ |\dot c^{(n)}(s_0)|}} $$
%
%and hence, recalling that $\ttt(s_0)=\dot c(s_0)$,
%$$\nnn(s_0)= \bbb(s_0)\tim\ttt(s_0)=\frac{c^{(n)}(s_0) }{ | c^{(n)}(s_0)|}\,. $$
%
\par In general, the binormal and the normal fail to be continuous at inflection points, see Example~\ref{Eflex}.
However, according to Proposition~\ref{PFS}, it turns out that they are both continuous when seen as functions in $\RP$.
This feature confirms that the natural ambient of definition for both the binormal and principal normal is indeed the projective plane $\RP$.
\smallskip\par Finally, in Sec.~\ref{Sec:force}, we define a measure $\calT$, that we call {\em torsion force}, that is obtained by performing the tangential variation of the length of the tangent indicatrix $\ttt_c$ that we have built up in Proposition~\ref{Ptan}.
Our torsion force may be compared with the {\em curvature force} $\cal K$ introduced in \cite{CFKSW}, that comes into the play by computing the first variation of the length
of curves with finite total curvature. In fact, in the smooth case we have:
$$ \calK=\gk\,\nnn\,d\calL^1 \,,\qquad \calT =k_{\#}\bigl(\gtau\,\bbb\,d\calL^1\bigr)  $$
where in the second formula we have set $k(s):=\int_0^s\gk(\l)\,d\l$, the primitive of the curvature of the curve.
\par In general, the curvature force $\cal K$ is a finite measure when the curve $c$ has finite total curvature $\TC(c)$, i.e., when the tantrix $\ttt=\dot c$ is a function of bounded variation. The torsion force $\calT$, instead, is a finite measure when the arc-length derivative of the tantrix
$\ttt_c$ from Proposition~\ref{Ptan} is a function with bounded variation. We shall see that this condition is satisfied if the curve $c$ has finite complete torsion
$\CT(c)$ in the sense of Alexandrov-Reshetnyak \cite{AR}.
\section{Weak binormal and total torsion of polygonals}\label{Sec:poly}
In this section, we introduce a weak notion of binormal indicatrix $\bbb_P$ for a polygonal $P$ in $\gR^3$, Definition~\ref{Dbin}. It is a rectifiable curve in the projective plane $\RP$ whose length is equal to the {\em total absolute torsion} of $P$.
%A notion of normal indicatrix $\nnn_P$ is consequently obtained.
\par
Let $P$ be a polygonal curve in $\gR^3$ with consecutive vertices $v_i$, $i=0,\ldots,n$, where $n\geq 3$ and $P$ is not closed, i.e., $v_0\neq v_n$. Without loss
of generality, we assume that every oriented segment $\s_i:=[v_{i-1},v_{i}]$ has positive length $\calL(\s_i):=\Vert v_{i}-v_{i-1}\Vert $, for $i=1,\ldots,n$, and that two consecutive segments are never aligned, i.e., the vector product
$\s_{i}\tim \s_{i+1}\neq 0_{\gR^3}$ for each $i=1,\ldots,n-1$.
\br\label{Ralign} If $\s_{i}\tim \s_{i+1}= 0_{\gR^3}$, we replace $\s_{i+1}$ with the oriented segment $[v_{i},v_{j+1}]$, where $j$ is the first index greater than $i$ such that $\s_{j}\tim \s_{j+1}\neq 0_{\gR^3}$. If $\s_{j}\tim\s_{j+1}= 0_{\gR^3}$ for each $j>i$, we set $b_i=b_{i-1}$ in definition \eqref{binormal} below.
\er
\adl\par\noindent{\large\sc Binormal vectors and torsion.}
In the definition by Penna \cite{Pe}, the {\em discrete unit binormal} is the unit vector given at each interior vertex $v_{i}$ of $P$ by the formula:
\beq\label{binormal} %b_{i-1}:={\s_{i-1}\tim \s_{i} \over |\s_{i-1}\tim \s_{i}| }\,, \qquad
b_i:=\frac{\s_{i}\tim \s_{i+1} }{ \Vert\s_{i}\tim \s_{i+1}\Vert }\,, \qquad i=1,\ldots, n-1\,. \eeq
The {\em torsion of $P$} is a function $\gtau(\s_i)$ of the interior oriented segments $\s_i$ defined as follows.
Let $i=2,\ldots,n-1$.
If the three segments $\s_{i-1},\,\s_i,\,\s_{i+1}$ are coplanar, i.e., if the vector product $b_{i-1}\tim b_i=0_{\gR^3}$, one sets $\gtau(\s_i)=0$.
Otherwise, one sets
$$ \gtau(\s_i):=\frac{\t_i}{ \calL(\s_i)}$$
where $\t_i$ denotes the angle between $-\p/2$ and $\p/2$ whose magnitude is the undirected angle between the binormals $b_{i-1}$ and $b_i$,
and whose sign is equal to the sign of the scalar product between the linearly independent vectors $b_{i-1}\tim b_i$ and $\s_i$.
Penna then defined the total torsion of $P$ through the sum:
$$ \sum_{i=2}^{n-1}\gtau(\s_i)\cdot\calL(\s_i)=\sum_{i=2}^{n-1}\t_i\,. $$
In a similar way, we define the {\em total absolute torsion} of $P$ by:
$$
\TAT(P):=\sum_{i=2}^{n-1}|\gtau(\s_i)|\cdot\calL(\s_i)=\sum_{i=2}^{n-1}|\t_i|\,. $$
\br\label{Rproj} In the above definitions, one considers angles between unoriented osculating planes. In fact, it may happen that the planes ${\text{span\,}}(\s_{i-1},\,\s_i)$ and ${\text{span\,}}(\s_i,\,\s_{i+1})$ are almost parallel,
but the directed angle between the binormal vectors $b_i$ and $b_{i+1}$ is equal to $\p-\e$ for some small $\e>0$. However, one gets $|\t_i|=\e$.
In facts, denoting by $\bullet$ the scalar product, in general one obtains
\beq\label{theta}
|\t_i|=\min \{\arccos(b_{i-1}\bullet b_{i}),\arccos(-b_{i-1}\bullet b_{i})\}\in[0,\p/2]\,. \eeq
\er
%\adl\par\noindent
{\large\sc An equivalent definition.}
In the classical approach by \cite{AR,Mi}, one considers the {\em tangent indicatrix} of $P$, i.e., the polygonal $\ttt_P$ in the Gauss sphere $\Sph^2$ obtained by letting $t_i:=\s_i/\calL(\s_i)\in\Sph^2$,
for $i=1,\ldots,n$, and connecting with oriented geodesic arcs $\g_i$ the consecutive points $t_{i}$ and $t_{i+1}$, for $i=1,\ldots,n-1$.
Therefore, one has $\calL(\g_i)=d_{\Sph^2}(t_i,t_{i+1})$, where $d_{\Sph^2}$ denotes the
geodesic distance on $\SP$.
\br\label{RTCP} The {\em total curvature $\TC(P)$ of $P$} is the sum of the {\em turning angles} $\a_i$ at the interior vertices of $P$, compare e.g. \cite{Su_curv}, and it is therefore equal to the length of $\ttt_P$, i.e.,
$$ \TC(P)=\sum_{i=1}^{n-1}\calL(\g_i)=\calL_{\SP}(\ttt_P)\,. $$
In particular, the arc-length parameterization $\ttt_P:[0,C]\to\Sph^2$, where $C:=\calL(\ttt_P)=\TC(P)$, is Lipschitz-continuous and piecewise smooth,
with $|\dot\ttt_P|=1$ everywhere except to a finite number of points, the edges of the tangent indicatrix $\ttt_P$,
which correspond to the interior segments of the polygonal $P$. \er
\br\label{Requal} With the previous assumptions on $P$, see Remark~\ref{Ralign}, the total absolute torsion of $P$ can be equivalently defined through the formula:
%Milnor then defined the total absolute torsion of $P$ through the formula:
%
$$\TAT(P):=\sum_{i=2}^{n-1}\wid\t_i $$
where $\wid\t_i\in[0,\p/2]$ is the shortest angle in $\Sph^2$ between the unoriented geodesic arcs $\g_{i-1}$ and $\g_i$ meeting at the edge $t_i$ of $\ttt_P$.
\par In fact, the geodesic arcs $\g_{i}$ are unique, as the consecutive points $t_i$ and $t_{i+1}$ are not antipodal. Moreover, we have $\wid \t_i=0$ exactly when $b_{i-1}\tim b_i= 0_{\gR^3}$, i.e., when $b_{i-1}=b_i$ or $b_{i-1}=-b_i$, so that $\gtau(\s_i)=0$.
Otherwise, we now check that $\wid \t_i=|\t_i|$ for each $i=1,\ldots,n-1$. By similarity, and up to a rotation, we can assume that $\s_i=(1,0,0)$. Setting
$\s_{i-1}=(\a_1,\be_1,\d_1)$ and $\s_{i+1}=(\a_2,\be_2,\d_2)$, one has $\s_{i-1}\tim\s_i=(0,\d_1,-\be_1)$ and $\s_{i}\tim\s_{i+1}=(0,-\d_2,\be_2)$, so that
$$ b_{i-1}=\frac{(0,\d_1,-\be_1)}{ \sqrt{\be_1^2+\d_1^2}}\,,\qquad b_i=\frac{(0,-\d_2,\be_2)}{ \sqrt{\be_2^2+\d_2^2}} $$
where $\s_{i-1},\,\s_i,\,\s_{i+1}$  are not coplanar provided that
%$\s_{i-1}\tim\s_i\neq 0_{\gR^3}$, $\s_{i}\tim\s_{i+1}\neq 0_{\gR^3}$, and
$b_{i-1}\tim b_i\neq 0_{\gR^3}$.
Now, the shortest angle $\wid\t_i$ between the geodesic arcs $\g_{i-1}$ and $\g_i$ meeting at $t_i$ is equal to the angle between the planes $\p_i^-$ and $\p_i^+$
spanned by the vectors $(\s_{i-1},\s_i)$ and $(\s_{i},\s_{i+1})$, respectively.
But the corresponding unit normals are $b_{i-1}$ and $b_i$, whence $\wid\t_i=|\t_i|$, where $|\t_i|$ is given by \eqref{theta}, as required. \er
\br\label{Rmeas} In an analytical approach, it turns out that the total curvature and absolute torsion of a polygonal $P$ can be seen as the total variation of mutually singular Radon measures $\m_P$ and $\n_P$ in $\gR^3$. In fact, with the above notation we have:
$$ \TC(P)=|\m_P|(\gR^3)\,,  \qquad \TAT(P)=|\n_P|(\gR^3) $$
where
$$ \m_P:=\sum_{i=1}^{n-1}\a_i\,\d_{v_i}\,, \qquad \ds \n_P:=\sum_{i=2}^{n-1}\t_i\,\Ha^1\pri{\s_i} $$
$\d_{v_i}$ being the unit Dirac mass at the vertex $v_i$ and $\Ha^1\pri\s_i$ the restriction to the segment $\s_i$ of the 1-dimensional Hausdorff measure $\Ha^1$. \er
\br\label{Rclose} If the polygonal $P$ is closed, i.e., $v_0=v_n$, the above notation is modified in a straightforward way: the torsion is defined
at all the $n$ segments $\s_i$, whereas the tangent indicatrix $\ttt_P$ is a closed polygonal curve in $\Sph^2$, so that $n$ angles are to be considered in both the equivalent definitions of $\TAT(P)$.
\er
{\large\sc The projective plane.} We have seen that the torsion is computed in terms of angles between undirected unit normal vectors $b_i$ of $\gR^3$, see Remarks~\ref{Rproj} and \ref{Requal}.
This implies that any reasonable notion of binormal (for non-smooth curves) naturally lives in the {\em real projective plane} $\RP$.
\par For this purpose, we recall that $\RP$ is defined by the quotient space
$\RP:=\Sph^2/\sim$, the equivalence relation being
$y\sim \wid y\iff y=\wid y$ or $y=-\wid y$, and hence the elements of $\RP$ are denoted by $[y]$.
The projective plane $\RP$ is naturally equipped with the induced metric
$$d_{\RP}([y],[\wid y]):=\min\{d_{\Sph^2}(y,\wid y),d_{\Sph^2}(y,-\wid y)\}\,. $$
Similarly to $(\SP,d_{\SP})$, the metric space
$(\RP,d_{\RP})$ is complete, and the projection map
$\Pi:\SP\to\RP$ such that $\Pi(y):=[y]$ is continuous.
Let $u:A\to\RP$ be a continuous map defined on an open set
$A\sb\Rn$. If $A\sb\Rn$ is simply connected, by the {\em
lifting theorem}, see e.g. \cite[p.~34]{Spa}, there
are exactly two continuous functions $v_i:A\to\SP$ such that
$[v_i]:=\Pi\circ v_i=u$, for $i=1,2$, with $v_2(x)=-v_1(x)$ for
every $x\in A$.
\par The manifold $\RP$ is non-orientable. Moreover, %following \cite[Sec.~3.1]{PR},
the mapping $g:\SP\to\gR^6$
$$%\beq\label{f}
g(y_1,y_2,y_3)=\Bigl(\frac{\sqrt
2}2 {y_1}^2,\,\frac{\sqrt 2}2 {y_2}^2,\,\frac{\sqrt 2}2
{y_3}^2,\,y_1 y_2,\, y_2 y_3,\,y_3 y_1\Bigr) $$%\eeq
induces an embedding
$$\wid g:\RP\to\wRP\,,\quad\wRP:=g(\SP)\sb\gR^6\,,\quad\wid g([y]):=g(y)\,. $$
Notice that $\wRP$\, is a non-orientable, smooth, compact,
connected submanifold of $\gR^6$ without boundary, such that
$|z|=\sqrt 2/2$ for every $z\in\wRP$. Also, $g$ maps the
equator $\SP\cap\{y^3=0\}$ into a circle $C$ of radius
$1/2$, covered twice, with constant velocity equal to one. The
circle $C$ is a minimum length generator of the first homotopy
group $\p_1(\wRP)\simeq \gZ_2$.
We also have $ \Ha^2(\wRP)=2\pi$, where $\Ha^2$ is the two-dimensional Hausdorff measure, compare e.g. \cite[Prop.~2.3]{MuLC}.
Moreover, $g$ is an isometric embedding. If e.g. a map $u:A\to\wRP$ is given by $u=g\circ v$ for
some smooth map $v:A\to\SP$, we in fact have
$$ \Vert D_iu\Vert^2=\Vert v\Vert^2\cdot\Vert D_iv\Vert^2+(v\bullet D_iv)^2 $$
for each partial derivative $D_i$. Therefore, since $\Vert v\Vert =1$
and $2\,(v\bullet D_iv)=D_i\Vert v\Vert^2=0$ a.e. for every $i$, we infer
that $\Vert Du\Vert =\Vert Dv\Vert $.
\adl\par\noindent{\large\sc Polar curve.} Using the above notation, and following Fenchel's approach \cite{Fe}, we now introduce the {\em polar} of the tangent indicatrix $\ttt_P$, a curve supported in the projective plane $\RP$,
in such a way that {\em the length in $\RP$ of the polar is equal to the total absolute torsion $\TAT(P)$}.
\par For this purpose, we recall that the support of $\ttt_P$ is the union of $n-1$ geodesic arcs $\g_i$, where $\g_i$ has initial point $t_i$ and end point
$t_{i+1}$, for $i=1,\ldots,n-1$. Since we assumed that consecutive segments of $P$ are never aligned, each arc $\g_i$ is non-trivial and well-defined.
According to the definition \eqref{binormal}, it turns out that the discrete unit binormal
$b_i\in\Sph^2$ is the ``north pole" corresponding to the great circle passing through $\g_i$
and with the same orientation as $\g_i$.
%For a point $b\in\Sph^2$, we denote $[b]\in\RP$ the corresponding equivalence class, and recall that $d_{\RP}([b_1],[b_2])=\min\{d_{\Sph^2}(b_1,b_2),\,d_{\Sph^2}(-b_1,b_2)\}$.
%
\par
\begin{figure}
	\centering
	\includegraphics[width=1.00\textwidth]{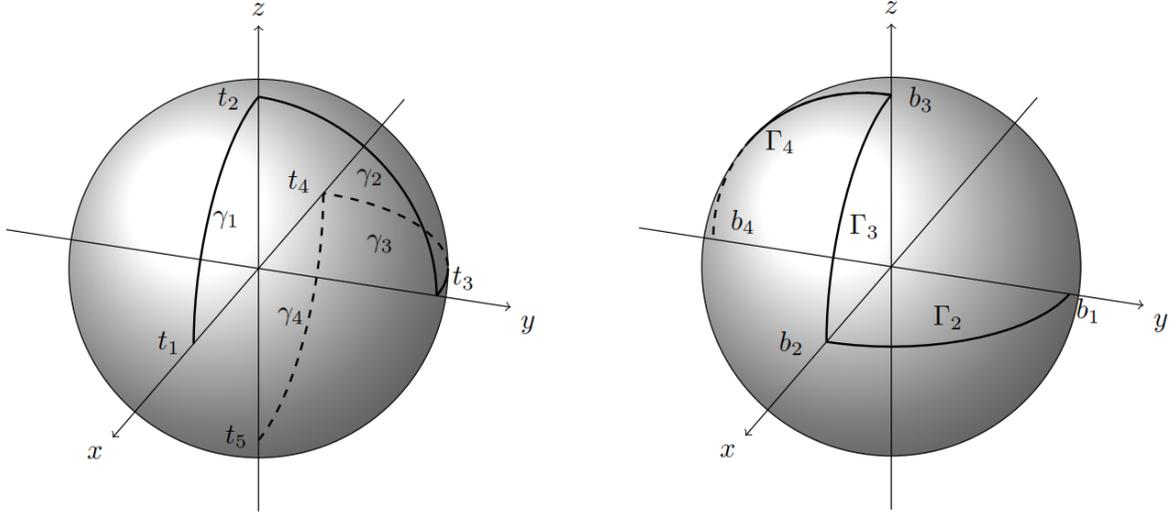}
		\caption{An example of a polygonal curve with tangent indicatrix moving as in the left figure. The weak binormal indicatrix moves as in the right figure. Since the weak binormal indicatrix lives in the projective space $\mathbb R\mathbb P^2$, in the figure we have drawn one of its two possible liftings into the sphere $\mathbb S^2$.}
\label{fig:SferaTB}
\end{figure}

 For any $i=2,\ldots,n-1$, we denote by $\GG_i$ the geodesic arc in $\RP$ with initial point $[b_{i-1}]$ and end point $[b_{i}]$.
Then $\GG_i$ is degenerate when $b_{i-1}=\pm b_i$, i.e., when the three segments $\s_{i-1},\,\s_i,\,\s_{i+1}$ are coplanar.
We thus have $\calL_{\RP}(\GG_i)=\wid\t_i=|\t_i|$ for each $i$, and hence that
$$\sum_{i=2}^{n-1}\calL_{\RP}(\GG_i)=\TAT(P)\,. $$
Also, for $i<n-2$ the end point of $\GG_i$ is equal to the initial point of $\GG_{i+1}$.
Finally, if $\TAT(P)=0$, i.e., if the polygonal $P$ is coplanar, all the arcs $\GG_i$ degenerate to a point $[b]\in\RP$, which actually identifies the binormal to $P$.
\bdf\label{Dpolar} {\em Polar} of the tangent indicatrix $\ttt_P$ is the oriented curve in $\RP$ obtained by connecting the consecutive geodesic arcs $\GG_i$,
for $i=2,\ldots,n-1$. \edf
{\large\sc Weak binormal.} Therefore, the polar of $\ttt_P$ connects by geodesic arcs in $\RP$ the consecutive discrete binormals $[b_i]$ of the polygonal $P$,
and its length is equal to the total absolute torsion $\TAT(P)$ of $P$.
In particular, it is a rectifiable curve.
This property allows us to introduce a suitable weak notion of binormal.
\bdf\label{Dbin} We denote {\em binormal indicatrix} of the polygonal $P$ the arc-length parameterization $\bbb_P$ of the polar in $\RP$ of the tangent indicatrix $\ttt_P$ (see Figure~\ref{fig:SferaTB}). \edf

\par We thus have $\bbb_P:[0,T]\to\RP$, where $T:=\calL_{\RP}(\bbb_P)=\TAT(P)$. Moreover, $\bbb_P$ is Lipschitz-continuous and piecewise smooth, with
$|\dot\bbb_P|=1$ everywhere except to a finite number of points.
\br\label{Rmonot} Differently from what happens for the length and the total curvature, the monotonicity formula fails to hold. More precisely, if $P'$ is a polygonal inscribed in $P$, by the triangular inequality we have $\calL(P')\leq \calL(P)$ and $\TC(P')\leq \TC(P)$, but it may happen that $\TAT(P')>\TAT(P)$.
This is due to the fact that the total absolute torsion of a polygonal $P$ can be computed as the sum of $\min\{\t_i,\pi-\t_i\}$, where $\t_i$ is the turning angle of the tantrix $\ttt_P$ at the $i$-th vertex.
\er
\bex\label{Emon} Let $P$ be a polygonal made of six segments $\s_i$, for $i=1,\ldots,6$, where the first three ones and the last three ones lay on two different planes $\Pi_1$ and $\Pi_2$.
Then the tantrix $\ttt_P$ connects with geodesic arcs in $\SP$ the consecutive points $t_i:=\s_i/\calL(\s_i)$, for $i=1,\ldots, 6$, where the triplets
$t_1,t_2,t_3$ and $t_4,t_5,t_6$ lay on two geodesic arcs, which are inscribed in the great circles corresponding to the vector spaces spanning the planes $\Pi_1$ and $\Pi_2$, respectively.
If both the angles $\a$ and $\be$ of $\ttt_P$ at the points $t_3$ and $t_4$ are small, then $\TAT(P)=\a+\be$.
\par Let $P'$ be the inscribed polygonal obtained by connecting the first point of $\s_3$ with the last point of $\s_4$.
The tantrix $\ttt_{P'}$ connects with geodesic arcs the consecutive points $t_1,t_2,w,t_5,t_6$, where the point $w$ lays in the minimal geodesic arc between $t_3$ and $t_4$.
Now, assume that the turning angle $\e$ of $\ttt_{P'}$ at the point $t_5$ satisfies $\a<\e<\p/2$, and that the two geodesic triangles with vertices $t_2,t_3,w$ and $w,t_4,t_5$ have the same area. By suitably choosing the position of the involved vertices, and by using the Gauss-Bonnet theorem in the computation, it turns out that $\TAT(P')-\TAT(P)=2(\e-\a)>0$.
\eex
%
%An important monotonicity property holds true. If $P$ and $P'$ are two polygonal curves in $\gR^3$,
%and $P'$ is obtained by replacing a segment $\s$ of $P$ with the two segments joining the end points of $\s$ with a new vertex, then:
%$$\calL(P)\leq \calL(P')\,,\qquad \TC(P)\leq\TC(P')\,,\qquad \TAT(P)\leq\TAT(P')\,. $$
%
%The first inequality is trivial. Moreover, looking at the tangent indicatrix and weak binormal corresponding to the polygonals, see Definition~\ref{Dbin},
%the triangle inequality in $\Sph^2$ and in $\RP$, respectively, yields that their lengths satisfy the inequalities:
%
%$$\calL_{\Sph^2}(\ttt_P)\leq \calL_{\Sph^2}(\ttt_{P'})\,, \qquad \calL_{\RP}(\bbb_P)\leq \calL_{\RP}(\bbb_{P'})\,. $$
%
\br\label{Rdual} For future use, we finally check the following inequality:
$$\TC_{\RP}(\bbb_P)\leq\calL_{\SP}(\ttt_{P})=\TC(P)\,. $$
In fact, for closed polygonals in the Gauss sphere such that three consecutive vertices never lie on the same geodesic, it turns out that polarity is an involutive transformation.
Therefore, the polar of (a lifting of) the binormal indicatrix $\bbb_P$ agrees with the polygonal in $\SP$ obtained by replacing
any chain of consecutive geodesic segments $\g_i$ of $\ttt_P$ which lay on some maximum circle, with a single geodesic arc obtained by connecting the end points of the chain.
In particular, the total curvature of $\bbb_P$ in $\RP$ is bounded by the length of $\ttt_P$.
\er
\section{Curves with finite total absolute torsion}\label{Sec:TAT}
In this section, we collect some notation concerning the total absolute torsion of curves in $\gR^3$.
We thus let $c$ be a curve in $\gR^3$ parameterized by
$c:I\to\gR^3$, where $I:=[a,b]$.
% and $c$ is continuous and one-to-one.
%
\par Any polygonal curve $P$
inscribed in $c$, say $P \ll c$, is obtained by choosing a
finite partition $\D:=\{a=\l_0<\l_1<\ldots<\l_{n-1}<\l_{n}=b\}$ of
\,$I$, say $P=P(\D)$, and letting $P:I\to\gR^3$ such that
$P(\l_i)=v_i:=c(\l_i)$ for $i=0,\ldots,n$, and $P(\l)$ is affine on
each interval $I_i:=[\l_{i-1},\l_i]$ of the partition, so that $P(I_i)=\s_i=[v_{i-1},v_i]$.
The {\em mesh} of the polygonal is defined by $\mesh P:=\sup\{\calL(\s_i)\mid i=1,\ldots,n \}$.
\par The length $\calL(c)$ and the total curvature $\TC(c)$
are respectively defined through the formulas:
$$\ba{rl} \calL(c):= &\sup\{\calL(P)\mid P\ll c\} \\ \TC(c):= &\sup\{\TC(P)\mid P\ll c\}\,. \ea $$
%$$  \TAT(c):= \sup\{\TAT(P)\mid P\ll c\}\,.  $$
%
\par Let $c$ be a curve in $\gR^3$ with {\em finite total curvature}, i.e., $\TC(c)<\i$. Then it is rectifiable, too, see e.g. \cite{Su_curv}.
Assume that $c:[0,L]\to\gR^3$ is its arc-length parameterization, whence $L=\calL(c)<\i$.
Since $c$ is a Lipschitz-continuous function, by Rademacher's theorem it is differentiable a.e. in $[0,L]$.
\par As a consequence, the tangent indicatrix $\ttt:[0,L]\to\Sph^2$ is well defined by setting $\ttt(s):=\dot c(s)$ for a.e. $s\in[0,L]$.
It is well-known that $\ttt$ is a function with bounded variation (see \cite{MuSa} for the notation on BV functions) and moreover that its essential variation in $\Sph^2$ agrees with the total curvature of $c$,
i.e., $\Var_{\Sph^2}(\ttt)=\TC(c)$. Notice that $\ttt$ is not continuous, as can be seen by taking a piecewise $C^1$ curve: a discontinuity point of $\ttt$ appears at any
edge point of $c$.
\par Moreover, by taking any sequence $\{P_h\}$ of inscribed polygonal curves such that
$\mesh P_h\to 0$, on account of Remark~\ref{Rmonot}, and by using a continuity argument, compare \cite{Su_curv}, one infers that $\calL(P_h)\to\calL(c)$ and $\TC(P_h)\to\TC(c)$.
%
%Notice that it suffices to take any sequence $\{\D_h\}$ of partitions of $I$ such that
%$\mesh \D_h\to 0$ and set $P_h=P(\D_h)$, since the uniform continuity of $c$ yields that $\mesh P_h\to 0$.
%
\smallskip\par\noindent{\large\sc Total absolute torsion.} Due to the lack of monotonicity described in Example~\ref{Emon}, we define the {\em total absolute torsion} $\TAT(c)$ of $c$ by means of the approach due to Alexandrov-Reshetnyak \cite{AR}.
\par
For this purpose, we recall that the {\em modulus} $\m_c(P)$ of a polygonal $P$ inscribed in $c$ is the maximum of the diameter of the arcs of $c$ determined by the consecutive vertices in $P$.
\par We also notice that if $c$ is a polygonal curve itself, there exists $\e>0$ such that any polygonal $P$ inscribed in $c$ and with modulus $\m_c(P)<\e$ satisfies
$\ttt_P=\ttt_c$, whence $\bbb_P=\bbb_c$ and definitely we get $\TAT(P)=\TAT(c)$.
It suffices indeed to take $\e$ lower than half of the mesh of the polygonal $c$, so that in every segment of $c$ there are at least two vertices of $P$.
\smallskip\par The above facts motivate the following definition:
\beq\label{DTAT}
 \TAT(c):=\lim_{\e\to 0^+}\sup\{ \TAT(P)\mid P\ll c \,,\,\,\m_c(P)<\e\}\,. \eeq
\par Therefore, if $\TAT(c)<\i$, for any sequence $\{P_h\}$ of polygonal curves inscribed in $c$ and satisfying $\m_c(P_h)\to 0$, one has $\sup_h\TAT(P_h)<\i$,
and one can find an optimal sequence as above in such a way that $\TAT(P_h)\to\TAT(c)$.
\smallskip\par Let now $c$ be a curve with finite total curvature and total absolute torsion.
%
%satisfying $\TAT(c)<\i$.
%Arguing as in \cite[Thm.~8.2.1]{AR}, it turns out that for any $s_0\in]0,L[$ such that any right-hand semi-neighborhood of $c(s_0)$ is not rectilinear,
%it is well defined a right-hand osculating plane at $c(s_0)$ in a suitable strong sense.
%This property implies that for each $\e>0$ we can find $\d>0$ such that for each polygonal $\wid P$ inscribed in the restriction $c_{\vert[s_0,s_0+\d]}$ of $c$
%to the interval $[s_0,s_0+\d]$, then $\TAT(\wid P)<\e$.
%
%As a consequence, on account of Remark~\ref{Rmonot}, and by using the same continuity argument as the one in \cite{Su_curv} for the total curvature result,
%one similarly obtains the following
%
In the next section, we shall see that it is possible to give a suitable weak notion of binormal indicatrix, a curve $\bbb_c$ in $\RP$
such that its length agrees with the total absolute torsion $\TAT(c)$, see \eqref{Varb} below.
\par As a consequence of Theorem~\ref{Tbin}, see Remark~\ref{Rappr}, we also obtain:
\bp\label{PTAT} Let $c$ be a curve in $\gR^3$ with both finite total curvature $\TC(c)$ and total absolute torsion $\TAT(c)$.
Then for any sequence $\{P_h\}$ of inscribed polygonal curves such that
$\m_c(P_h)\to 0$, one has $\TAT(P_h)\to\TAT(c)$.
\ep
\par For this purpose, we first discuss here the regular case, i.e., when curvature and torsion are defined as in the usual way.
\adl\par\noindent{\large\sc The smooth case.} Let $c$ be a smooth regular curve in $\gR^3$ defined through the arc-length parameterization (so that $|\dot c|=1$ a.e.).
Assuming $\ddot c\neq 0$ everywhere, and letting $\ttt:=\dot c$, $\nnn:=\dot\ttt/|\dot\ttt|$, $\gk:=|\dot\ttt|$, $\bbb:=\ttt\tim\nnn$,
the classical Frenet-Serret formulas for the spherical frame $(\ttt,\nnn,\bbb)$ of $c$ give:
\beq\label{FS} \dot\ttt=\gk\,\nnn\,,\qquad \dot\nnn=-\gk\,\ttt+\gtau\,\bbb\,,\qquad\dot\bbb=-\gtau\,\nnn \eeq
where $\gk$ is the (positive) curvature and $\gtau$ the torsion of the curve.
\par By Proposition~\ref{PTAT}, and on account of the density result from \cite[Prop.~4]{Pe}, one readily obtains:
\bc\label{CTAT} If $c$ is a smooth regular curve in $\gR^3$,
then
$$ \TAT(c)=\int_0^L|\gtau(s)|\,ds\,. $$
\ec
\br\label{RFrenet} Notice that a rectifiable curve may have unbounded total curvature but zero torsion (just consider a planar curve).
Conversely, by taking $s\in[0,1]$ and letting $\gk(s)\equiv 1$ and $\gtau(s)=(1-s)^{-1}$, solutions to the Frenet-Serret system \eqref{FS} are rectifiable curves $c$ such that $\int_c \gk\,ds=1$ but $\int_c|\gtau|\,ds=+\i$.
\er
\par As the following example shows, the (absolute value of the) torsion may be seen as the curvature of the tantrix, when computed in the sense of the spherical geometry.
\bex\label{Ehelix} Given $R>0$ and $K\geq 0$, we let $c:[-L/2,L/2]\to\gR^3$ denote the helicoidal curve
$$c(s):=(R\cos(s/v),R\sin(s/v),K s/(2\p v))\,, \quad  s\in[-L/2,L/2] $$
where we denote $v:=(R^2+(K/2\p)^2)^{1/2}$ and choose $L:=2\pi v$, so that $|\dot c|\equiv 1$ and the length $\calL(c)=L$. Moreover,
$c(\pm L/2)=(\pm R,0,\pm K/2)$, and $c(0)=(R,0,0)$.
We thus have
$$ \ba{rl} \ttt(s)= & v^{-1}(-R\sin(s/v),R\cos(s/v),K/2\p) \\
\nnn(s)= & (-\cos(s/v),-\sin(s/v),0) \\
\bbb(s)= &v^{-1}((K/2\p)\sin(s/v),-(K/2\p)\cos(s/v),R) \ea $$
so that both curvature and torsion are constant, $\gk\equiv Rv^{-2}$, $\gtau\equiv v^{-2}(K/2\p)$. Therefore,
the integral of the curvature and of the torsion of $c$ are readily obtained:
$$\int_c \gk\,ds=L\cdot\gk=\frac{2\p R}{ v}\,, \qquad  \int_c |\gtau|\,ds=L\cdot\gtau=\frac{K}{ v}\,, \qquad v:=(R^2+(K/2\p)^2)^{1/2}. $$

We now compute the {\em spherical curvature} $\gk_{\Sph^2}(\ttt)$ of the tantrix $\ttt$, a closed curve embedded in the Gauss sphere $\Sph^2$ and parameterizing (when $K>0$) a small circle whose radius depends on $R$ and $K$. We consider a sequence of (strongly converging) polygonal curves $\{\ttt_n\}$ in $\Sph^2$ inscribed in the tantrix $\ttt$. The total curvature of $\ttt_n$ is equal to the sum of the width in $\Sph^2$ of the angles between consecutive segments. When $n\to\i$, by uniform convergence we obtain the total curvature of $\ttt$ in $\Sph^2$. Actually, it agrees with the integral of the absolute torsion of $c$, i.e.,
$$ \int_\ttt  \gk_{\Sph^2}(\ttt)\,ds =\frac{K}{ v}= \int_c |\gtau|\,ds\,.$$

To this purpose, for each $n\in\Nat^+$, we let $t_n(i):=\ttt(s_i)$, where $s_i=(L/n)i$ and $i\in\gZ\cap[-n,n]$, and we consider the closed spherical polygonal generated by the consecutive points $t_n(i)\in \Sph^2$.

The turning angle in $\Sph^2$ of two consecutive geodesic segments $t_n(i-1)t_n(i)$ and $t_n(i)t_n(i+1)$, agrees with the angle between the two planes in $\gR^3$ spanned by $0_{\gR^3}$ and the end points of the above segments, i.e., between the normals $t_n(i-1)\tim t_n(i)$ and $t_n(i)\tim t_n(i+1)$.
By symmetry, such an angle $\t_n$ does not depend on the choice of $i$, and will be computed at $i=0$.
The total spherical curvature of the polygonal being equal to $n\cdot \t_n$, we check:
$$ \lim_{n\to \i} n\cdot \t_n = \frac{K}{ v}\,. $$
\par In fact, in correspondence to the middle point we have
$$t_n(0)=v^{-1}(0,R,K/2\p)\,, \quad t_n(\pm 1)=v^{-1}(\mp R\sin(2\p/n),R\cos(2\p/n),K/2\p) $$
so that we get
$$ t_n(0)\tim t_n(\pm 1)=\frac{R}{ v^{2}}\cdot\Bigl({ K\over2\p}\,\Bigl(1-\cos{2\p\over n}\Bigr), \mp { K\over2\p}\,\sin{2\p\over n},
\pm R\,\sin{2\p\over n}  \Bigr)\,.$$
Denoting for simplicity
$$ M_n:=\Vert t_n(0)\tim t_n(\pm 1)\Vert =\frac{R}{ v^{2}}\cdot\bigl( (K/2\p)^2 2(1-\cos(2\p/n))+R^2\sin^2(2\p/n)  \bigr)^{1/2} $$
and setting $N^\pm_n:=\pm (t_n(0)\tim t_n(\pm 1))/M_n$, %$N^-_n:=P_n(-1)\tim P_n(0)/|P_n(-1)\tim P_n(0)|$,
we compute
$$ N^+_n \tim N^-_n =\frac{R^2}{ {M_n}^{2}}\,(K/2\p)\,\sin(2\p/n)\,2(1-\cos(2\p/n))\cdot(0,-R,(K/2\p))$$
$$ \Vert N^+_n \tim N^-_n\Vert =\frac{R^2}{ {M_n}^{2}}\,(K/2\p)\,\sin(2\p/n)\,2(1-\cos(2\p/n))\,v\,.$$
By symmetry, the turning angle of the geodesic arcs connecting two consecutive points $t_n(i)$ does not depend on the choice of $i$ and is equal to $$\t_n:=\arcsin\Vert N^+_n \tim N^-_n\Vert \,. $$

Since for $n\to\i$ we have $2(1-\cos(2\p/n))\sim (2\p/n)^2$ and $\sin(2\p/n)\sim 2\p/n$, we get
$M_n\sim R (2\p/n) v$ and finally $n\cdot\t_n\sim n\,\Vert N^+_n \tim N^-_n\Vert \to K/v$ where, we recall, $\int_c|\gtau|\,ds =K/v$.
\eex
\br In the previous example, we have considered a sequence $\{\ttt_n\}$ of polygonal curves in $\Sph^2$ inscribed in the tantrix $\ttt$ of $c$ and converging to $\ttt$ in the sense of the Hausdorff distance.
In general, each $\ttt_n$ is not the tangent indicatrix of a polygonal inscribed in $c$.
However, the total spherical curvature $n\cdot\t_n$ of $\ttt_n$ clearly agrees with the length in $\RP$ of the polar of $\ttt_n$, which is constructed as in
Sec.~\ref{Sec:poly}, see Definition~\ref{Dpolar}.
\par Now, one may similarly consider a sequence $\{P_h\}$ of polygonals inscribed in $c$, each one made of $h$ segments with the same length, so that $\mesh P_h\to 0$.
The total absolute torsion $\TAT(P_h)$ of $P_h$ agrees with the
length in $\RP$ of the binormal indicatrix $\bbb_{P_h}$, see Definition~\ref{Dbin}. By means of a similar computation (that we shall omit), one can show that
$\calL_{\RP}(\bbb_{P_h})\to K/v$ as $h\to\i$, in accordance with the formula in Corollary~\ref{CTAT}.
\er
\section{Weak binormal of a non-smooth curve}\label{Sec:bin}
In this section, we consider rectifiable curves $c$ in $\gR^3$ with finite total curvature $\TC(c)$ and finite (and non zero) total absolute torsion $\TAT(c)$.
Using a density approach by polygonals, in Theorem~\ref{Tbin} we show that a {\em weak} notion of {\em binormal indicatrix} of $c$ is well-defined.
For smooth curves, we shall recover the classical binormal, see Theorem~\ref{Tbinsmooth} and Remark~\ref{Rbin}.
Finally, similar properties concerning the tangent indicatrix are discussed in Propositions~\ref{Ptan} and \ref{Ptansmooth}.
\par More precisely, we shall define a Lipschitz-continuous function $\bbb_c:[0,T]\to\RP$, where $T=\TAT(c)$, satisfying $|\dot\bbb_c|=1$ a.e. in $[0,T]$.
Therefore, $\bbb_c$ is a curve in $\RP$ with length equal to the total absolute torsion of $c$, i.e.,
\beq\label{Varb}\calL_{\RP}(\bbb_c)=\TAT(c)\,. \eeq
\par This is the content of our first main result:
\bt\label{Tbin} Let $c$ be a curve in $\gR^3$ with finite total curvature $\TC(c)$ and finite (and non-zero) total absolute torsion $T:=\TAT(c)$. There exists a rectifiable curve $\bbb_c:[0,T]\to\RP$ parameterized by arc-length, so that $\calL_{\RP}(\bbb_c)=\TAT(c)$, satisfying the following property. For any sequence $\{P_h\}$ of inscribed polygonal curves, let $b_h:[0,T]\to\RP$ denote for each $h$ the parameterization with constant velocity of the
binormal indicatrix $\bbb_{P_h}$ of $P_h$, see Definition~$\ref{Dbin}$.
If $\m_c(P_h)\to 0$, then $b_h\to \bbb_c$ uniformly on $[0,T]$ and $\calL_{\RP}(b_h)\to\calL_{\RP}(\bbb_c)$. \et
\br\label{Rappr} Recalling that $\calL_{\RP}(b_h)=\TAT(P_h)$, Proposition~\ref{PTAT} readily follows.
\er
\par Furthermore, we shall see that if $c$ is smooth in the sense of the previous section (so that the Frenet-Serret formulas \eqref{FS} hold),
the binormal $\bbb(s)$ of $c$ agrees with the value of a suitable lifting of the weak binormal $\bbb_c$ in $\Sph^2$, when computed at the expected point.
\bt\label{Tbinsmooth} Let $c:[0,L]\to\gR^3$ be a rectifiable curve of class $C^3$ parameterized in arc-length, so that $L=\calL(c)$.
Assume that $\ddot c(s)\neq 0$ for each $s\in[0,L]$, so that the spherical frame $(\ttt,\nnn,\bbb)$ of $c$ is well-defined.
Let $\bbb_c:[0,T]\to\RP$ be the rectifiable curve in $\RP$ defined in Theorem~$\ref{Tbin}$, so that $T=\TAT(c)$.
Then, for each $s\in]0,L[$ there exists $t(s)\in[0,T]$ such that
$$\bbb(s)=\wid\bbb_c(t(s)) $$
for a unique lifting $\wid\bbb_c$ of $\bbb_c$ in $\SP$.
Moreover, $t(s)$ is equal to the total absolute torsion $\TAT(c_{\vert [0,s]})$ of the curve $c_{\vert [0,s]}:[0,s]\to\gR^3$. In particular, we have:
%
%In addition, if the torsion $\gtau$ of $c$ never vanishes on $[0,L]$, then $t=t(s)$ is determined by the bijective transition function $t:[0,L]\to[0,T]$
%
\beq\label{torsion} t(s)=\int_0^s|\gtau(\l)|\,d\l \qquad\fa\,s\in[0,L] \eeq
where $\gtau(\l)$ is the torsion of the curve $c$ at the point $c(\l)$. \et
\br\label{Rbin} Notice that if the torsion $\gtau$ of $c$ (almost) never vanishes, the function $t(s):[0,L]\to[0,T]$ in equation \eqref{torsion} is strictly increasing, and its inverse $s(t):[0,T]\to[0,L]$ gives
$$ \wid\bbb_c(t)=\bbb(s(t))\qquad \fa\,t\in[0,T]\,,\qquad T=\TAT(c)\,.  $$
Therefore, in this case, the weak binormal $\bbb_c$ in $\RP$, when suitably lifted to $\SP$, agrees with the {\em arc-length parameterization of the binormal $\bbb$} of $c$. \er
{\large\sc Tangent indicatrix.} Similarly to Theorems~\ref{Tbin} and \ref{Tbinsmooth}, we also obtain the following
properties concerning the tantrix.
\bp\label{Ptan} Let $c$ be a curve in $\gR^3$ with finite total curvature $C:=\TC(c)$ and with no points of return.
% and finite total absolute torsion $\TAT(c)$.
Then, there exists a rectifiable curve $\ttt_c:[0,C]\to\SP$, parameterized by arc-length, so that $\calL_{\SP}(\ttt_c)=\TC(c)$, satisfying the following property. For any sequence $\{P_h\}$ of inscribed polygonal curves such that $\mesh P_h\to 0$, denoting by $t_h:[0,C]\to\SP$ the parameterization with constant velocity of the
tangent indicatrix $\ttt_{P_h}$ of $P_h$, then $t_h\to \ttt_c$ uniformly on $[0,C]$ and $\calL_{\SP}(t_h)\to\calL_{\SP}(\ttt_c)$. \ep
\br\label{Runique} If $c$ has points of return, i.e., if e.g. for some $s\in]0,L[$ we have $\ttt(s-)=-\ttt(s+)$, then the curve $\ttt_c$ is uniquely determined up to the choice of
the geodesic arc in $\SP$ connecting $\ttt(s-)$ and $\ttt(s+)$.
\er
\bp\label{Ptansmooth} Let $c:[0,L]\to\gR^3$ be a curve of class $C^2$ parameterized in arc-length, so that $L=\calL(c)$, and let
$\ttt_c:[0,C]\to\SP$ be the rectifiable curve in $\SP$ defined in Proposition~$\ref{Ptan}$, so that $C=\TC(c)$.
Then, for each $s\in]0,L[$ there exists $k(s)\in[0,C]$ such that the tangent indicatrix $\ttt:=\dot c$ satisfies
$$\ttt(s)=\ttt_c(k(s)) \,. $$
Moreover, $k(s)$ is equal to the total curvature $\TC(c_{\vert [0,s]})$ of the curve $c_{\vert [0,s]}:[0,s]\to\gR^3$, whence:
\beq\label{curvature} k(s)=\int_0^s\gk(\l)\,d\l \qquad\fa\,s\in[0,L] \eeq
where $\gk(\l):=\Vert \ddot c(\l)\Vert $ is the curvature of $c$ at the point $c(\l)$. \ep
\br As before, if the curvature $\gk$ of $c$ (almost) never vanishes, the function $k(s):[0,L]\to[0,C]$ in equation \eqref{curvature} is strictly increasing, and its inverse $s(k):[0,C]\to[0,L]$ gives
$$ \ttt_c(k)=\ttt(s(k))\qquad \fa\,k\in[0,C]\,,\qquad C=\TC(c)\,.  $$
%
%Whence, the weak tangent $\ttt_c$ agrees with the arc-length parameterization of the tantrix $\ttt$.
\er
{\large\sc Proofs.} We now give the proofs of the previous results.
%
%\adl\par\noindent
\bpff{\sc of Theorem~\ref{Tbin}}: It is divided into four steps.
\smallskip\par\noindent{\sc Step~1.} Choose an optimal sequence $\{P_h\}$ of polygonal curves inscribed in $c$ such that $\m_c(P_h)\to 0$ and
$T_h\to T$, where $T_h:=\TAT(P_h)$ and $T=\TAT(P)$.
For $h$ large enough so that $T_h>0$, the binormal indicatrix of $P_h$ has been defined by the arc-length parameterization
$\bbb_{P_h}:[0,T_h]\to\RP$ of the curve in $\RP$ given by the polar of the tangent indicatrix $\ttt_{P_h}$, see Definition~\ref{Dbin}. Whence it is a rectifiable curve such that
$\calL_{\RP}(\bbb_{P_h})=T_h$ and $\Vert \dot\bbb_{P_h}\Vert =1$ a.e. on $[0,T_h]$.
%
%Since $\mesh P_h\to 0$, by Proposition~\ref{PTAT} we also know that $T_h\to T^-$, where $T:=\TAT(c)$.
%
\par Define $b_h:[0,T]\to\RP$ by $b_h(s):=\bbb_{P_h}((T_h/T)s)$,
so that $\Vert \dot
b_h(s)\Vert = T_h/T$ a.e., where $T_h/T\to 1$. By Ascoli-Arzela's
theorem, we can find a subsequence $\{b_{h_k}\}$ that uniformly
converges in $[0,T]$ to some Lipschitz continuous function $b:[0,T]\to\RP$, and we denote $b=\bbb_c$.
\smallskip\par\noindent{\sc Step~2.} We claim that $\dot b_h\to \dot b=\dot\bbb_c$ strongly in $L^1$.
As a consequence, we deduce that $\Vert \dot \bbb_c\Vert =1$ a.e. on $[0,T]$, and hence that
$$ \calL_{\RP}(\bbb_c)=\int_0^T\Vert \dot\bbb_c(s)\Vert \,ds=T=\TAT(c)\,. $$
%
%that equality holds true in the above formula. In particular, we obtain that $=\TAT(c)$, as required.
%
\par In order to prove the claim, recalling from Sec.~\ref{Sec:poly} that $\wid g:\RP\to\wRP\sb\gR^6$
is the isometric embedding of the projective plane, we shall denote here $\ul f:=\wid g\circ f$, for any function $f$ with values in $\RP$, and we consider
the tantrix $\tt_{h}$ of the curve ${\ul {b_h}}:[0,T]\to\wRP$, i.e., $\tt_{h}(s)=\dot{\ul {b_h}}(s)/\Vert \dot{\ul {b_h}}(s)\Vert $.
We have $\calL_{\RP}(b_h)=\TAT(P_h)$ and $\Vert \dot{\ul {b_h}}(s)\Vert =T_h/T$, whereas by Remark~\ref{Rdual}
$$\TC_{\RP}(b_h)\leq\calL_{\SP}(\ttt_{P_h})=\TC(P_h)\,. $$
Therefore, it turns out that the essential total variation of $\tt_h$ in $\wRP$ is lower than the sum $\TC(P_h)+\TAT(P_h)$. We thus get:
$$\sup_h\Var_{\wRP}(\tt_h)\leq \TC(c)+\TAT(c)<\i\,. $$
As a consequence, by compactness, a subsequence of $\{\dot{\ul {b_h}}\}$ converges weakly-* in the $\BV$-sense to some $\BV$-function $v:[0,T]\to\wRP$.
\par We show that $v(s)=\dot{\ul {b}}(s)$ for a.e. $s\in[0,T]$. This yields that the sequence $\{\dot b_{h}\}$ converges strongly in $L^1$ (and hence a.e. on $[0,T]$) to the function $\dot b$.
\par In fact, using that by
Lipschitz-continuity
$$ {\ul {b_h}}(s)={\ul {b_h}}(0)+\int_0^s \dot{\ul {b_h}}(\l)\,d\l\qquad \fa\, s\in [0,T]
$$
and setting
$$ V(s):={\ul {b}}(0)+\int_0^s v(\l)\,d\l\,,\qquad s\in [0,T] $$
by the weak-* $\BV$ convergence $\dot{\ul {b_h}}\wc v$, which implies the
strong $L^1$ convergence, we have ${\ul {b_h}}\to V$ in $L^\i$,
hence $\dot{\ul {b_h}}\to \dot V=v$ a.e. on $[0,T]$. But we already
know that ${\ul {b_h}}\to {\ul {b}}$ in $L^\i$, thus we get $v=\dot{\ul {b}}$.
\smallskip\par\noindent{\sc Step~3.} Let now $\{\wid P_h\}$ denote any sequence of polygonal curves inscribed in $c$ such that $\m_c(\wid P_h)\to 0$. We claim that
possibly passing to a subsequence, the binormals $\bbb_{\wid P_h}$ converge uniformly (up to reparameterizations) to the curve $\bbb_c$.
\par In fact, we recall that the polar of the tantrix $\ttt_{P}$ to a polygonal curve $P$ is defined in terms of vector products of couples of consecutive points of its geodesic segments,
the vector product being continuous. Moreover, the {\em Frech\'et distance} (see e.g. \cite[Sec.~1]{Su_curv}) between the two sequences $\{\ttt_{P_h}\}$ and $\{\ttt_{\wid P_h}\}$ goes to zero.
This property follows from the equiboundedness of the total curvatures.
Whence, the polars of $\ttt_{P_h}$ and of $\ttt_{\wid P_h}$ must converge uniformly (up to reparameterizations) to the same limit function.
Therefore, the sequence $\bbb_{\wid P_h}$ converges in the Frech\'et distance to the curve $\bbb_c$ obtained in Step~1.
\smallskip\par\noindent{\sc Step~4.}
Now, if $\{\wid P_h\}$ is the (not relabeled) subsequence obtained in Step~3, by repeating the argument in Step~1 we infer that
the limit function $b=\bbb_c$ is unique. As a consequence, a contradiction argument yields that all the sequence $\{b_{h}\}$ uniformly converges to $\bbb_c$
and that the limit curve $\bbb_c$ does not depend on the choice of the sequence $\{P_h\}$ of
inscribed polygonals satisfying $\m_c(P_h)\to 0$. Therefore, the curve $\bbb_c$ is
identified by $c$. Arguing as in Step~2, we finally infer that $\calL_{\RP}(b_h)\to\calL_{\RP}(\bbb_c)$, as required.
\epff
\bpff{\sc of Theorem~\ref{Tbinsmooth}}: For any given $s\in]0,L[$, since $\Vert \dot c(s)\Vert =1$ and $\ddot c(s)\neq 0$, the binormal is defined by $\bbb(s):=\ttt(s)\tim\nnn(s)$, with $\ttt(s):=\dot c(s)$ and
$\nnn(s):=\ddot c(s)/\Vert \ddot c(s)\Vert $,
so that $\dot c(s)\tim\ddot c(s)\neq 0$ and
$$ \bbb(s)=\frac{\dot c(s)\tim\ddot c(s) }{ \Vert \dot c(s)\tim\ddot c(s) \Vert  }\,. $$
We thus may and do choose a sequence of polygonals $\{P_h\}$ inscribed in $c$ such that $\m_c(P_h)\to 0$ and (with the notation from Sec.~\ref{Sec:poly} for $P=P_h$)
the following properties hold for any $h\in\Nat^+$ large enough : \ben
\item the four points $v_{i-2}=c(s-2h)$, $v_{i-1}=c(s-h)$, $v_i=c(s+h)$, $v_{i+1}=c(s+2h)$ are consecutive (and interior) vertices of $P_h$;
\item the three segments $\s_{i-1}=v_{i-1}-v_{i-2}$, $\s_i=v_{i}-v_{i-1}$, $\s_{i+1}=v_{i+1}-v_{i}$ satisfy $\s_{i-1}\tim\s_i\neq 0_{\gR^3}$ and $\s_{i}\tim\s_{i+1}\neq 0_{\gR^3}$.
\een
\par By taking the second order expansions of $c$ at $s$, we get
$$ \ba{rl} \s_{i-1}= & \ds -\dot c(s)\,h+\frac{3}{ 2}\,\ddot c(s)\,h^2+o(h^2)\,, \\
\s_{i}= & \ds {2}\,\ddot c(s)\,h^2+o(h^2)\,, \\
\s_{i+1}= & \ds \dot c(s)\,h+\frac{3}{ 2}\,\ddot c(s)\,h^2+o(h^2) \ea
$$
and hence
$$  \s_{i-1}\tim\s_i = 2h^2\,\ddot c(s)\tim \dot c(s)+o(h^3)\,, \qquad \s_{i}\tim\s_{i+1} = 2h^2\,\ddot c(s)\tim \dot c(s)+o(h^3)\,.$$
On account of \eqref{binormal}, we thus get for any $h$ large:
$$b_{i-1}(h)=\frac{\s_{i-1}\tim \s_{i} }{ \Vert \s_{i-1}\tim \s_{i}\Vert  }=-\bbb(s)+o(h^3),\,\,\,
b_i(h)=\frac{\s_{i}\tim \s_{i+1} }{ \Vert \s_{i}\tim \s_{i+1}\Vert  }=-\bbb(s)+o(h^3) $$
so that in particular $b_i(h)\to -\bbb(s)$ as $h\to\i$.
\par Now, consider the polygonal $P_h(s)$ given by the union of the segments $\s_1,\ldots,\s_{i-1},\s_i$ of $P_h$. It turns out that the total absolute torsion of $P_h(s)$ satisfies $\TAT(P_h(s))=t_h(s)$ for some number $t_h(s)\in[0,\TAT(P_h)]$. Since $\TAT(P_h)\to\TAT(c)\in\gR^+$, possibly passing to a subsequence, the sequence $\{t_h(s)\}$ converges to some number $t(s)\in[0,T]$.
By Theorem~\ref{Tbin}, we thus infer that $b_i(h)\to \bbb_c(t(s))$ as $h\to\i$, whence we obtain $\bbb(s)=-\bbb_c(t(s))$.
\par Moreover, since both the end points of the segment $\s_i$ of $P_h$ converge to $c(s)$ as $h\to\i$, whereas $\m_c(P_h)\to 0$, by Proposition~\ref{PTAT} we deduce that
$\TAT(P_h(s))\to\TAT(c_{\vert [0,s]})$, which yields the equality $t(s)=\TAT(c_{\vert [0,s]})$.
Since by smoothness of the curve $c$
$$ \TAT(c_{\vert [0,s]})=\int_0^s\Vert \dot\bbb(\l)\Vert \,d\l$$
recalling that $\dot \bbb(\l)=-\gtau(\l)\,\nnn(\l)$, we finally obtain the equality \eqref{torsion}.
\epff
\bpff{\sc of Proposition~\ref{Ptan}}: Following the proof of Theorem~\ref{Tbin}, we choose $h$ large enough so that $C_h:=\TC(P_h)>0$,
and we denote by $\ttt_{P_h}:[0,C_h]\to\SP$ the arc-length parameterization of the tantrix $\ttt_{P_h}$, so that
$C_h=\calL_{\SP}(\ttt_{P_h})$ and $\Vert \dot\ttt_{P_h}\Vert =1$ a.e. on $[0,C_h]$. Since $\mesh P_h\to 0$, we have $C_h\to C^-$, where $C:=\TC(c)$.
Setting $t_h:[0,C]\to\SP$ by $t_h(s):=\ttt_{P_h}((C_h/C)s)$, as in Step~1 we can find a subsequence $\{t_{h_k}\}$ that uniformly
converges in $[0,C]$ to some Lipschitz continuous function $t:[0,C]\to\SP$.
Moreover, as in Steps~3-4 we deduce that $t$ does not depend on the choice of $\{P_h\}$, and that all the sequence $\{t_{h}\}$ uniformly converges to $t$, so that the curve $\ttt_c:=t$ is identified by $c$.
\par We claim that $\calL_\SP(\ttt_c)=C$. As a consequence, since the equality $\Vert \dot t_h\Vert =C_h/C$ a.e. yields that $\Vert \dot \ttt_c\Vert \leq 1$ a.e., whereas
$\calL_\SP(\ttt_c)=\int_0^C \Vert \dot \ttt_c(s)\Vert \,ds$, we infer that $\Vert \dot \ttt_c\Vert =1$ a.e., as required.
\par It remains to prove the claim. Since $\ttt=\dot c$ is a function of bounded variation, for each $h$ we can find a partition $\D_h$ of $[0,L]$ in $2^h$ intervals $I^h_i=[s^h_{i-1},s^h_i]$, for $i=1,\ldots,2^h$, satisfying the following properties:
\ben \item $\D_{h+1}$ is a refinement of $\D_h$, and $\mesh \D_h\to 0$ as $h\to\i$\,;
\item for each $i$, the end points of the intervals $I^h_i$ are Lebesgue points of $\ttt$, with Lebesgue values $\ttt(s^h_{i-1})$ and $\ttt(s^h_i)$\,;
\item if $f_h:[0,L]\to\SP$ is the piecewise constant function with $f_h(s)=\ttt(s^h_i)$ for each $s\in]s^h_{i-1},s^h_i[$ and each $i$, then $\Var_{\gR^3}(f_h)\to\Var_{\gR^3}(\ttt)$\,.
\een
\par Let now $\g_h$ denote the spherical polygonal in $\SP$ obtained by connecting the consecutive vertices $\ttt(s^h_i)$. Then, $\calL_\SP(\g_h)=\Var_\SP(\g_h)\to\Var_\SP(\ttt)=\TC(c)$. On the other hand, the Frech\'et distance between the two sequences $\{\ttt_{P_h}\}$ and $\{\g_h\}$ goes to zero.
Therefore, $\g_h$ converges to $\ttt_c$ in the Frech\'et distance. As a consequence, each polygonal $\g_h$ is inscribed in $\ttt_c$, which yields that
$\calL_\SP(\g_h)\to\calL_\SP(\ttt_c)$, and hence that $\calL_\SP(\ttt_c)=\TC(c)$, which completes the proof. \epff
%
%\bpf{\em of Proposition~\ref{Ptan}} As in the proof of Theorem~\ref{Tbin}, but with $\ttt_{P_h}$, $C_h$, $\SP$, $C$, $t_h$, $t$, and $\ttt_c$ instead of
%$\bbb_{P_h}$, $T_h$, $\RP$, $T$, $b_h$, $b$, and $\bbb_c$, respectively, this time using that
%the essential total variation in $\SP$ of the tantrix $\tt_h$ of $t_h$ is lower than the complete torsion $\CT(P_h)$ in the sense of \cite{AR}, which is lower than the sum $\TC(P_h)+\TAT(P_h)$. We omit any further detail.
%\epf
%
\br\label{Rtan} It turns out that the essential total variation in $\SP$ of the tantrix $\tt_h$ of $t_h$ is lower than the complete torsion $\CT(P_h)$ in the sense of \cite{AR}.
Therefore, if in addition  the curve $c$ has finite complete torsion in the sense of \cite{AR}, $\CT(c)<\i$, as in Step~3 of the proof of Theorem~\ref{Tbin} we infer that the derivative $\dot\ttt_c$ is a function of bounded variation, and that $\dot t_h$ converges to $\dot\ttt_c$ weakly-* in the $\BV$-sense, and hence a.e. in $[0,C]$.
We finally recall that a curve with finite total curvature and total absolute torsion may have infinite complete torsion.
\er
\bpff{\sc of Proposition~\ref{Ptansmooth}}: Similarly to the proof of Theorem~\ref{Tbinsmooth}, for any $s\in]0,L[$ we choose $\{P_h\}$ inscribed in $c$ such that $\mesh P_h\to 0$ and for any $h\in\Nat^+$ the two points
$v_{i-1}=c(s-h)$ and $v_{i}=c(s+h)$ are consecutive (and interior) vertices of $P_h$.
We thus get $\s_i:=v_{i}-v_{i-1}=2\dot c(s)\,h+o(h)$, whence
$t_i(h):=\s_i/\Vert \s_i\Vert \to \dot c(s)=\ttt(s)$ as $h\to\i$.
Also, denoting again by $P_h(s)$ the polygonal corresponding to the segments $\s_1,\ldots,\s_{i-1},\s_i$ of $P_h$, we have $\TC(P_h(s))=k_h(s)\in[0,\TC(P_h)]$, where $\TC(P_h)\to C\in\gR^+_0$, whence a subsequence of
 $\{k_h(s)\}$ converges to some $k(s)\in[0,C]$.
Proposition~\ref{Ptan} yields that $t_i(h)\to \ttt_c(k(s))$ as $h\to\i$, whence we get $\ttt(s)=\ttt_c(k(s))$.
We clearly have $\TC(P_h(s))\to\TC(c_{\vert [0,s]})$, which implies that
$$k(s)=\TC(c_{\vert [0,s]})=\ds\int_0^s\Vert \dot\ttt(\l)\Vert \,d\l\,. $$
Recalling that $\dot \ttt=\gk\,\nnn$, we finally obtain the equality \eqref{curvature}. \epff
\section{Weak normal of a non-smooth curve}\label{Sec:nor}
We have seen that the curvature of an open polygonal $P$ is a non-negative measure $\m_P$ concentrated at the interior vertices of $P$, whereas
the torsion is a signed measure $\n_P$ concentrated at the interior segments, see Remark~\ref{Rmeas}.
Since these two measures are mutually singular, in principle there is no analogous to the classical formula by Fenchel for the (principal) normal of smooth curves in $\gR^3$, namely
\beq\label{normal} \nnn=\bbb\tim\ttt\,. \eeq
\par In this section, following Banchoff \cite{Ba}, a weak notion of normal indicatrix of a polygonal is introduced, Definition~\ref{Dnor}, in such a way that formula \eqref{normal} continues to hold.
As a consequence, according to the cited Fenchel's approach, the principal normal of a curve with finite total curvature and absolute torsion is well-defined in a weak sense, Theorem~\ref{Tnor}.
\adl\par\noindent{\large\sc Weak normal of polygonals.}
Let $P$ be an open polygonal in $\gR^3$ with non-degenerate segments.
Denoting $C=\TC(P)$ and $T=\TAT(P)$, we first choose two suitable curves
$$ \wid\ttt_P:[0,C+T]\to\RP\,, \qquad \wid\bbb_P:[0,C+T]\to\RP $$
which on one side inherit the properties of the tangent indicatrix and of the binormal indicatrix of $P$, respectively, and on the other side take account of the order in which curvature and torsion are defined along $P$.
More precisely, we shall recover the properties
\beq\label{Pnew}  \calL_{\RP}(\wid\bbb_P)=\TC_{\RP}(\wid\ttt_P)=\TAT(P)\,,\qquad \TC_{\RP}(\wid\bbb_P)\leq\calL_{\RP}(\wid\ttt_P)=\TC(P)
\eeq
(where all equalities hold in the case of closed polygonals),
which are satisfied (up to a lifting) by the curves $\ttt_P$ and $\bbb_P$ defined in Sec.~\ref{Sec:poly}.
Moreover, in accordance to the mutual singularities of the measures $\m_P$ and $\n_P$, see Remark~\ref{Rmeas},
one curve is constant when the other one parameterizes a geodesic arc, whose length is equal to the curvature or to the (absolute value of the) torsion at one
vertex or segment of $P$, respectively.
\par Recalling the notation from Sec.~\ref{Sec:poly}, we let $v_i$, $i=0,\ldots,n$, denote the vertices, and $\s_i:=[v_{i-1},v_{i}]$, $i=1,\ldots,n$, the oriented segments of $P$. Also, we let $t_i:=\s_i/\calL(\s_i)\in\Sph^2$,
for $i=1,\ldots,n$, and $\g_i$ is a minimal geodesic arc in $\SP$ connecting the consecutive points $t_{i}$ and $t_{i+1}$, for $i=1,\ldots,n-1$.
Notice that $\g_i$ is unique when $t_{i+1}\neq -t_i$, and it is trivial when $t_{i+1}= t_i$.
Finally, $\GG_i$ is the geodesic arc in $\RP$ with initial point $[b_{i-1}]$ and end point $[b_{i}]$, for any $i=2,\ldots,n-1$, where $b_i$ is the discrete binormal \eqref{binormal}. Therefore, $\GG_i$ is trivial when $b_i=\pm b_{i-1}$.
We thus have
$$ \TC(P)=\sum_{i=1}^{n-1}\calL_{\SP}(\g_i)\,,\qquad \TAT(P)=\sum_{i=2}^{n-1}\calL_{\RP}(\GG_i)\,.$$
\br\label{Rprod} In order to explain our construction below, let us choose a lifting $\wih\bbb_P:[0,T]\to\Sph^2$ of the (continuous) curve $\bbb_P$ from Definition~\ref{Dbin},
and let $\wih b_i$ and $\wih\GG_i$ denote the points and geodesic arcs corresponding to $[b_i]$ and $\GG_i$.
For $i=1,\ldots,n-1$, we let $\wid \g_i=\wih b_i\tim \g_i$, i.e., $\wid \g_i$ is the oriented geodesic arc in $\Sph^2$ obtained by means of the vector product of the lifted discrete binormal $\wih b_i$ with each point in the support of the arc $\g_i$. For $i=2,\ldots,n-1$, we also let $\wid \GG_i=\wih\GG_i\tim t_{i+1}$, i.e., $\wid \GG_i$ is the oriented geodesic arc in $\Sph^2$ obtained by means of the vector product of each point in the support of the lifted arc $\wih\GG_i$ with the direction $t_{i+1}$.
\par
It turns out that for $i=1,\ldots,n-2$, the final point of $\wid\g_i$ agrees with the initial point of $\wid\GG_{i+1}$, and that the final point of $\wid\GG_{i+1}$ agrees with the initial point of $\wid\g_{i+1}$.
Using this order to join the geodesic arcs, one obtains a rectifiable curve in $\Sph^2$ whose total length is equal to the sum of the lengths of $\ttt_P$
and of $\bbb_P$, i.e., to $\TC(P)+\TAT(P)$.
However, since the curve depends on the chosen lifting of the binormal, it is more natural to work in the projective plane. Therefore, we shall consider the geodesic arcs $[\g_i]:=\Pi(\g_i)$ with end points $[t_i]:=\Pi(t_i)$, where $\Pi:\SP\to\RP$ is the canonical projection.
\er
\par Recalling that $C:=\TC(P)$ and $T=\TAT(P)$, we shall denote for brevity $C_0:=0$, $T_1:=0$, and
$$C_i:=\sum_{j=1}^i\calL_{\SP}(\g_j)\,, \,\,\, i=1,\ldots,n-1\,,\quad T_i:=\sum_{j=2}^{i}\calL_{\RP}(\GG_j)\,,\,\,\, i=2,\ldots,n-1 \,.$$
Notice that $C_i=C_{i-1}$ if $\g_i$ is trivial, i.e., when $t_{i+1}=t_i$, and that $T_i=T_{i-1}$ when $\GG_i$ is trivial, i.e., when $b_i=\pm b_{i-1}$.
\par
We define $\wid \ttt_P:[0,C+T]\to\RP$ and $\wid \bbb_P:[0,C+T]\to\RP$ as follows:
\ben \item $\wid\ttt_P$ parameterizes with velocity one the oriented geodesic arc $[\g_i]$ on the interval $[C_{i-1}+T_i,C_i+T_i]$, for $i=1,\ldots, n-1$
such that $\g_i$ is non-trivial;
\item $\wid\ttt_P$ is constantly equal to $[t_i]$ on the interval $[C_{i-1}+T_{i-1},C_{i-1}+T_{i}]$, for $i=2,\ldots, n-2$\,;
\item $\wid\bbb_P$ is constantly equal to $[b_i]$ on the interval $[C_{i-1}+T_i,C_i+T_i]$, for $i=1,\ldots, n-1$\,;
\item  $\wid\bbb_P$ parameterizes with velocity one the oriented geodesic arc $\GG_i$ on the interval $[C_{i-1}+T_{i-1},C_{i-1}+T_{i}]$, for $i=2,\ldots, n-2$
such that $\GG_i$ is non-trivial.
\een
\par
\begin{figure}
\centering
	\includegraphics[width=0.60\textwidth]{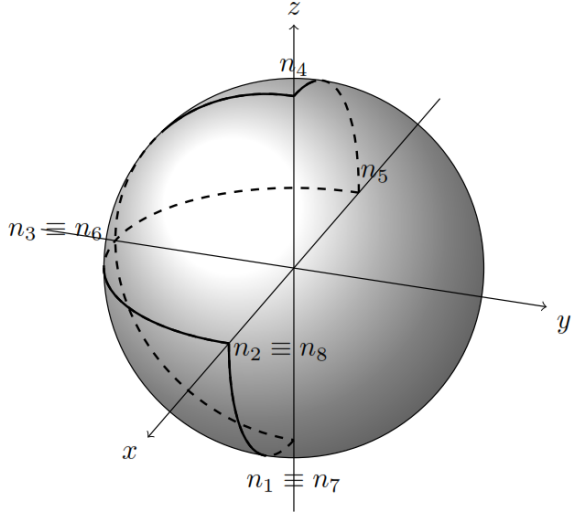}
	\caption{The weak normal indicatrix of the curve whose tangent and binormal indicatrix are those in Figure~\ref{fig:SferaTB} of page \pageref{fig:SferaTB}. Again, for the sake of the illustration we consider one of the two liftings of the normal indicatrix into the sphere ${\mathbb S}^2$.}
\label{fig:SferaN}
\end{figure}

 The functions $\wid\ttt_P$ and $\wid\bbb_P$ are both continuous, and property \eqref{Pnew} is readily checked.
%We thus may and do choose a lifting $\wih\bbb_P:[0,C+T]\to\SP$ of the continuous function $\wid\bbb_P$, so that $\wih\bbb_P$ is a curve in $\SP$ such that
%$\Pi\circ\wih\bbb_P=\wid\bbb_P$.
Furthermore, it turns out
that {\em the unit vectors $\wid\ttt_P(s)$ and $\wid\bbb_P(s)$ are orthogonal, for a.e. $s\in[0,C+T]$.}
As a consequence, we are able to define the weak normal according to the formula \eqref{normal}.
\bdf\label{Dnor} {\em Normal indicatrix} of the polygonal $P$ is the curve $\nnn_P:[0,C+T]\to\RP$ (see Figure~\ref{fig:SferaN}) given by the pointwise vector product
$$ \nnn_P(s):=\wid \bbb_P(s)\tim \wid\ttt_P(s)\in\RP\,, \qquad s\in[0,T+C]\,. $$ \edf
\par For closed polygonals, the above notation is modified in a straightforward way, arguing as in Remark~\ref{Rclose}.
\br By the definition, it turns out that
$$ \calL_{\RP}(\nnn_P)=\calL_{\RP}(\wid\ttt_P)+\calL_{\RP}(\wid\bbb_P)=\TC(P)+\TAT(P)\,. $$
Notice that, the curvature and torsion of $P$ being mutually singular measures, see Remark~\ref{Rmeas}, the above equality is the analogous in the category of polygonals to the integral formulas
$$ \ba{c} \ds \int_c \Vert \dot\nnn(s)\Vert \,ds=\int_c\sqrt{\gk^2(s)+\gtau^2(s)}\,ds\,, \\
\ds \int_c \gk(s)\,ds=\TC(c)\,,\qquad \int_c |\gtau(s)|\,ds=\TAT(c)\ea $$
for smooth curves $c$, which clearly follow from the Frenet-Serret formulas \eqref{FS}.
\par
Moreover, we have $\Vert \dot\nnn_P(s)\Vert =1$ for a.e. $s\in[0,C+T]$. In fact, by the definition of $\wid\ttt_P$ and $\wid\bbb_P$, we get:
\ben \item for $i=1,\ldots, n-1$ and $s\in]C_{i-1}+T_i,C_i+T_i[$, we have $\wid\bbb_P(s)\equiv [b_i]\in\RP$ and hence $\dot\nnn_P(s)= [b_i]\tim \dot{\wid\ttt}_P(s)$,
where $\Vert \dot{\wid\ttt}_P(s)\Vert =1$ and $[b_i]$ is orthogonal to $\dot{\wid\ttt}_P(s)$, if $\g_i$ is non-trivial;
\item for $i=2,\ldots, n-2$ and $s\in]C_{i-1}+T_{i-1},C_{i-1}+T_{i}[$ we have $\wid\ttt_c(s)\equiv [t_i]$ and hence
$\dot\nnn_P(s)= \dot{\wid\bbb}_P(s)\tim[t_i]$, where $\Vert \dot{\wid\bbb}_P(s)\Vert =1$ and $[t_i]$ is orthogonal to $\dot{\wid\bbb}_P(s)$, if $\GG_i$ is non-trivial. \een
\er
\br\label{Rnor} Notice that the turning angle in $\RP$ of the curve $\nnn_P$ is equal to $\p/2$ at each ``non-trivial" vertex of $\nnn_P$.
Indeed, from a vertex of $\nnn_P$ we move by rotating either around $\ttt_\alpha$ or $\bbb_\beta$ ($\beta=\alpha$ or $\beta=\alpha-1$), where $\ttt_\alpha\perp\bbb_\beta$, hence the two curves are orthogonal.
More precisely, for $i=1,\ldots, n-1$, if both the geodesic arcs $[\g_i]$ and $\GG_{i+1}$ are non-degenerate, they meet orthogonally at the vertex $\nnn_P(C_i+T_i)$ of $\nnn_P$.
Similarly, for any $i=2,\ldots, n-2$ such that both the geodesic arcs $\GG_{i+1}$ and $[\g_{i+1}]$ are non-degenerate, they meet orthogonally at the vertex
$\nnn_P(C_i+T_{i-1})$.
\er
%
%\adl\par\noindent
{\large\sc Weak normal of curves.} In the same spirit as in Theorem~\ref{Tbin}, for non-smooth curves (that may have points of return or planar pieces) we now obtain our second main result.
In view of Remark~\ref{Rtan}, we need the stronger assumption that the curve has finite complete torsion $\CT(c)$ in the sense of \cite{AR}.
To this purpose, we recall that the implication $\CT(c)<\i \Lra \TAT(c)<\i$ holds true in general, whereas the
implication $\CT(c)<\i \Lra \TC(c)<\i$ is satisfied provided that the curve has no points of return.
\bt\label{Tnor} Let $c$ be a curve in $\gR^3$ with finite total curvature $C:=\TC(c)$, finite complete torsion $\CT(c)$, and finite total absolute torsion $T:=\TAT(c)$. There exists a rectifiable curve $\nnn_c:[0,C+T]\to\RP$ parameterized by arc-length, so that $\calL_{\RP}(\nnn_c)=C+T$, satisfying the following property. For any sequence $\{P_h\}$ of inscribed polygonal curves, let $n_h:[0,C+T]\to\RP$ denote the parameterization with constant velocity of the
normal indicatrix $\nnn_{P_h}$ of $P_h$, see Definition~$\ref{Dnor}$.
If $\m_c(P_h)\to 0$, then $n_h\to \nnn_c$ uniformly on $[0,C+T]$ and $\calL_{\RP}(n_h)\to\calL_{\RP}(\nnn_c)$. \et
\bpf We clearly may and do assume that each $P_h$ has non-degenerate segments. By Definition~\ref{Dnor}, setting $C_h=\TC(P_h)$ and $T_h=\TAT(P_h)$, the normal indicatrix of $P_h$ is the curve $\nnn_{P_h}:[0,C_h+T_h]\to\RP$ given by
$ \nnn_{P_h}(s):=\wid \bbb_{P_h}(s)\tim \wid\ttt_{P_h}(s)$, so that
$ \calL_{\RP}(\nnn_{P_h})=C_h+T_h$, and $\Vert \dot \nnn_{P_h}\Vert =1$ a.e. on $[0,C_h+T_h]$.
Also, condition $\m_c(P_h)\to 0$ yields that $C_h\to C$ and $T_h\to T$.
\par Setting $n_h:[0,C+T]\to\RP$ by $n_h(s):=\nnn_{P_h}((C_h+T_h)s/(C+T))$,
as before we deduce that possibly passing to a subsequence, the sequence $\{n_{h}\}$ uniformly
converges to some Lipschitz continuous function $\nnn_c:[0,C+T]\to\RP$.
\par We claim that $\Vert \dot \nnn_c\Vert = 1$ a.e. in  $[0,C+T]$. This yields that
$$ \calL_{\RP}(\nnn_c)=\int_0^{C+T}\Vert \dot \nnn_c(s)\Vert \,ds= C+T=\TC(c)+\TAT(c)\,. $$
For this purpose, we note that by Definition~\ref{Dnor} we have $n_h(s)=\wid b_h(s)\tim\wid t_h(s)$ for each $s\in[0,T+C]$, where
$$\wid b_h(s):=\wid\bbb_{P_h}((C_h+T_h)s/(C+T))\,, \qquad \wid t_h(s):=\wid\ttt_{P_h}((C_h+T_h)s/(C+T)) \,. $$
As in Theorem~\ref{Tbin} and Proposition~\ref{Ptan}, using that (by Remark~\ref{Rdual}) we again have:
$$\calL_{\RP}(\wid b_h)=\TC_{\RP}(\wid t_h)=\TAT(P_h)\,,\qquad \TC_{\RP}(\wid b_h)\leq\calL_{\RP}(\wid t_h)=\TC(P_h)\,, $$
we deduce that (possibly passing again to a subsequence) $\wid b_h\to \wid b$ and $\wid t_h\to \wid t$ strongly in $L^1$ (and uniformly) to some continuous functions with bounded variation
$\wid b,\wid t:[0,C+T]\to\RP$, and that the approximate derivatives $\dot{\wid b_h}\to \dot{\wid b}$ and $\dot{\wid t_h}\to\dot{\wid t}$ a.e. on $[0,C+T]$, see Remark~\ref{Rtan}.
This yields that $\nnn_c(s)=\wid b(s)\tim\wid t(s)$ and hence:
$$ \ba{rl}\ds\lim_{h\to\i} \dot n_h(s)= & \ds\lim_{h\to\i}\bigr(\dot{\wid b}_h(s)\tim\wid t_h(s)+\wid b_h(s)\tim\dot{\wid t}_h(s) \bigr) \\
= & \ds
\bigl(\dot{\wid b}(s)\tim\wid t(s)+\wid b(s)\tim\dot{\wid t}(s)\bigr)=\dot\nnn_c(s) \ea $$
for a.e. $s\in[0,C+T]$. But we already know that $\Vert \dot n_h(s)\Vert =(C_h+T_h)/(C+T)$ for a.e. $s$, where
$C_h\to C$ and $T_h\to T$, whence the claim is proved.
\par We now show that the limit function $\nnn_c$ does not depend on the initial choice of the approximating sequence $\{P_h\}$.
As a consequence, as before we conclude that the weak normal $\nnn_c$ only depends on $c$, and that the whole sequence $\{n_{h}\}$ converges to $\nnn_c$.
\par In fact, if we choose another sequence of polygonals $\{P^{(1)}_h\}$, we know that the sequences $\{\wid b_{P^{(1)}_h}\}$ and
$\{\wid t_{P^{(1)}_h}\}$ take the same limit as the one of the sequences
$\{\wid b_{P_h}\}$ and $\{\wid t_{P_h}\}$, respectively.
Moreover, the corresponding limit function $\nnn_c^{(1)}$ has length equal to the length of $\nnn_c$ on each interval $I\sb[0,C+T]$, and hence the same property holds true for the corresponding couples of functions $\wid b$, $\wid b^{(1)}$ and $\wid t$, $\wid t^{(1)}$, respectively. These facts imply that $\nnn_c^{(1)}=\nnn_c$, as required.
\epf
\br On account of Remark~\ref{Rnor}, denoting by $\tt_h$ the tantrix of the curve $\ul n_h:=\wid g(n_h)$ in $\wRP$, in general we do not have $\sup_h \Var_{\wRP}(\tt_h)<\i$.
Therefore, we cannot argue as in Theorem~\ref{Tbin} to conclude that the sequence $\dot n_h$ converges weakly in the $\BV$-sense (and hence strongly in $L^1$) to the function $\dot \nnn_c$. Actually, the derivative $\dot\nnn_c$ of the weak normal $\nnn_c$ is not a function with bounded variation, in general.
 \er
{\large\sc The case of smooth curves.} We finally have:
\bp\label{Pnorsmooth} Let $c:[0,L]\to\gR^3$ be a smooth curve satisfying the hypotheses of Theorem~$\ref{Tbinsmooth}$, so that $L=\calL(c)$, $C=\TC(c)$, and $T=\TAT(c)$ are finite.
Let $s:[0,C+T]\to[0,L]$ be the inverse of the increasing and bijective function $t:[0,L]\to[0,C+T]$ given by
\beq\label{ts} t(s):=\int_0^s(\gk(\l)+|\gtau(\l)|)\,d\l\,, \qquad s\in[0,L] \eeq
where $\gk(\l)$ and $\gtau(\l)$ are the curvature and torsion of the curve $c$ at the point $c(\l)$.
Then the principal normal $\nnn$ in $\SP$ of the curve $c$, and the curve $\nnn_c$ in $\RP$ given by Theorem~$\ref{Tbin}$, are linked by the formula:
\beq\label{nst} [\nnn(s(t))]=\nnn_c(t)\in\RP\qquad \fa\,t\in[0,C+T]\,. \eeq
\ep
\bpf For any given $s\in]0,L[$, we choose a sequence $\{P_h\}$ as in the proof of Theorem~\ref{Tbinsmooth}, and we correspondingly denote:
$$ b_i(h):=\frac{\s_{i}\tim \s_{i+1} }{ \Vert \s_{i}\tim \s_{i+1}\Vert  }\,, \qquad t_i(h):=\frac{\s_i }{\Vert \s_i\Vert }\,. $$
Letting $t_h(s):=\TC(P_h(s))+\TAT(P_h(s))$, this time we infer that (possibly passing to a subsequence)
$t_h(s)\to t(s):=\TC(c_{\vert [0,s]})+\TAT(c_{\vert [0,s]})$, so that $t(s)$ satisfies the formula \eqref{ts}.
As a consequence, arguing as in the proofs of Theorem~\ref{Tbinsmooth} and Proposition~\ref{Ptansmooth}, on account of Theorem~\ref{Tnor} this time we get:
$$ \lim_{h\to\i}[b_i(h)\tim t_i(h)]= \nnn_c(t(s))\,. $$
Since $b_i(h)\to \bbb(s)$ and $t_i(h)\to\ttt(s)$, we also have $b_i(h)\tim t_i(h)\to \nnn(s)$, so that formula \eqref{nst} holds. We omit any further detail.
\epf
\section{On the spherical indicatrices of smooth curves}\label{Sec:smooth}
The trihedral $(\ttt,\nnn,\bbb)$ is well-defined everywhere in the case of regular curves $\g$ in $\gR^3$ of class $C^2$ such that $\g'(t)$ and $\g''(t)$ are always independent vectors,
and the Frenet-Serret formulas \eqref{FS} hold true if in addition $\g$ is of class $C^3$.
\par Fenchel in \cite{Fe} used a geometric approach in order to define (under weaker hypotheses on the curve) the osculating plane. He chooses the binormal $\bbb$ as a smooth function.
Therefore, the principal normal is the smooth function given by $\nnn=\bbb\tim\ttt$. The Frenet-Serret formulas continue to hold, but this time the curvature may vanish and even be negative.
He also calls $\gk$-inflection or $\gtau$-inflection a point of the curve where the curvature or the torsion changes sign, respectively.
\par By using an analytical approach, we recover some of the ideas by Fenchel in order to define the binormal (and principal normal).
In general, it turns out that the binormal and normal fail to be continuous at the inflection points (see Example~\ref{Eflex}).
However, both the binormal and normal are continuous when seen as functions in the projective plane $\RP$.
\par For this purpose, in the sequel we shall assume that $\g:[a,b]\to\gR^3$ satisfies the following properties:
\ben \item $\g$ is differentiable at each $t\in[a,b]$ and $\g'(t)\neq 0_{\gR^3}$, i.e., $\g$ is a regular curve;
\item for each $t_0\in]a,b[$, the function $\g$ is of class $C^n$ in a neighborhood of $t_0$, for some $n\geq 2$, and
$\g^{(n)}(t_0)\neq 0_{\gR^3}$, but $\g^{(k)}(t_0)=0_{\gR^3}$ for $2\leq k\leq n-1$, if $n\geq 3$.
\een
\par We thus denote by $c(s):=\g(t(s))$ the arc-length parameterization of the curve $\g$, i.e., $t(s)=s(t)^{-1}$, with
$s(t):=\int_a^t\Vert \dot\g(\l)\Vert \,d\l\in[0,L]$, where $L:=\calL(\g)$.
\bp\label{PFS} Under the above assumptions, the Frenet-Serret frame $(\ttt,\bbb,\nnn)$ is well-defined for each $s_0\in[0,L]$ by:
\beq\label{tbn} \ba{l} \ds \ttt(s_0):=\dot c(s_0)\,, \qquad \bbb(s_0):=\frac{\dot c(s_0)\tim c^{(n)}(s_0) }{ \Vert c^{(n)}(s_0)\Vert }\,,  \\
 \ds \nnn(s_0):=\bbb(s_0)\tim\ttt(s_0)=\frac{ c^{(n)}(s_0) }{ \Vert c^{(n)}(s_0)\Vert } \ea \eeq
where $s_0=s(t_0)$ and $n\geq 2$ is given as above. Furthermore, $\ddot c(s_0)=0_{\gR^3}$ at a finite or countable set of points, and if
$\ddot c(s_0)\neq 0_{\gR^3}$, then $\nnn(s_0)=\ddot c(s_0)/\Vert \ddot c(s_0)\Vert $.
Finally, $[\bbb]$ and $[\nnn]$ are continuous functions with values in $\RP$.\ep
\bpf We set $\ttt(s):=\dot c(s)$ for each $s$. If $\ddot c(s_0)=0_{\gR^3}$,
then for some $n\geq 3$ and for $h$ small (and non-zero) we have
\beq\label{eqn-der}\ddot c(s_0+h)=c^{(n)}(s_0)\,\frac{h^{n-2}}{(n-2)!}+o(h^{n-2})\,.\eeq
This implies that $\ddot c(s)=0_{\gR^3}$ only at isolated points $s\in[0,L]$, hence at an at most countable set.
\par If $\ddot c(s_0)\neq 0_{\gR^3}$, one defines as usual $\nnn(s_0):=\ddot c(s_0)/\Vert \ddot c(s_0)\Vert $ and $\bbb(s_0):=\dot c(s_0)\tim\ddot c(s_0)/\Vert \ddot c(s_0)\Vert $.
In fact, the orthogonality property $\dot c(s_0)\bullet \ddot c(s_0)=0$ yields that $(\dot c(s_0)\tim \ddot c(s_0))\tim \dot c(s_0)=\ddot c(s_0)$.
\par If $\ddot c(s_0)=0_{\gR^3}$, for $h$ small we have $\dot c(s_0+h)\bullet \ddot c(s_0+h)=0_{\gR^3}$. Letting $h\to 0$, we obtain
that
$\dot c(s_0)\bullet c^{(n)}(s_0)=0_{\gR^3}$, whence $\dot c(s_0)\tim c^{(n)}(s_0)\neq 0_{\gR^3}$ and $\Vert \dot c(s_0)\tim c^{(n)}(s_0)\Vert =\Vert c^{(n)}(s_0)\Vert >0$.
As a consequence, the binormal is well-defined at $s_0$ such that $\ddot c(s_0)=0_{\gR^3}$ by the limit
$$ \bbb(s_0)=\lim_{h\to 0}\bbb(s_0+h)=\lim_{h\to 0}\frac{\dot c(s_0+h)\tim \ddot c(s_0+h) }{ \Vert \dot c(s_0+h)\tim \ddot c(s_0+h)\Vert }=
\frac{\dot c(s_0)\tim c^{(n)}(s_0) }{ \Vert c^{(n)}(s_0)\Vert }$$
and the principal normal is then defined by letting $ \nnn(s_0):=\bbb(s_0)\tim\ttt(s_0)$, where this time the orthogonality property $\dot c(s_0)\bullet c^{(n)}(s_0)=0_{\gR^3}$ yields that
$$ \nnn(s_0)= \frac{(\dot c(s_0)\tim c^{(n)}(s_0))\tim \dot c(s_0) }{ \Vert c^{(n)}(s_0)\Vert }=\frac{c^{(n)}(s_0) }{ \Vert c^{(n)}(s_0)\Vert }\,. $$
\par Finally, we observe that where $\ddot c\neq0_{\gR^3}$, both $\nnn$ and $\bbb$ are continuous (as functions valued in $\mathbb S^2$, hence also as functions valued in $\RP$). Therefore, the problematic points are where $\ddot c=0_{\gR^3}$, which is a set of isolated points. At one of these point, $\nnn(s_0)$ is ideally given by the limit of $\ddot c(s_0+h)/\Vert \ddot c(s_0+h)\Vert $, as $h\to0$. Using equation \eqref{eqn-der}, it is easy to see that, depending on the parity of the derivative order $n$, either the right and left limits coincide (thus the limit exists, and $\nnn$ is continuous at $s_0$) or they are opposite to one another.
Hence $\nnn$ and $\bbb$ may not be continuous as sphere-valued functions, but they are continuous as projective-valued function, since their directions are well defined and continuous. \epf
\br If in addition we assume that $\g$ is of class $C^3$, it turns out that the Frenet-Serret formulas \eqref{FS} hold true outside the at most countable set of inflection points. In fact, we have seen that $\ddot c(s)=0_{\gR^3}$ only at isolated points $s\in[0,L]$.
\er
\bex\label{Eflex} Let $c:[-1,1]\to\gR^3$ be a regular curve with derivative
$$\dot c(s)={1\over\sqrt 2}\,\bigl(1,s^2,\sqrt{1-s^4}\bigr),\quad s\in[-1,1] $$
so that $\Vert\dot c(s)\Vert\equiv 1$ and hence $\ttt(s)=\dot c(s)$. For $s\in]-1,1[$, we compute
$$\ddot c(s)={\sqrt 2 s\over \sqrt{1-s^4}}\,\bigl(0,\sqrt{1-s^4},-s^2\bigr)\,,\quad c^{(3)}(s)=\sqrt 2\,\Bigl(0,1,{s^2(s^4-3)\over (1-s^4)^{3/2}} \Bigr)\,. $$
Therefore, if $0<|s|<1$ we have $\ddot c(s)\neq 0_{\gR^3}$ and hence
$$ \nnn(s)={s\over |s|}\,\bigl(0,\sqrt{1-s^4},-s^2\bigr)\,,\quad \bbb(s)={s\over |s|}\,{1\over\sqrt 2}\,\bigl(-1,s^2,\sqrt{1-s^4}\bigr)\,. $$
In particular, the normal and binormal can be extended by continuity at $s=\pm 1$ by letting
$\nnn(\pm 1):=(0,0,\mp 1)$ and $\bbb(\pm1):=2^{-1/2}(\mp 1,\pm 1,0)$.
\par
Furthermore, for $0<|s|<1$ we get:
$$ \gk(s):=\Vert \ddot c(s)\Vert={\sqrt 2|s|\over \sqrt{1-s^4}}\,,\quad\gtau(s):={\bigl(\dot c(s)\tim\ddot c(s)\bigr)\bullet c^{(3)}(s) \over \Vert \ddot c(s)\Vert^2}=-{\sqrt 2s\over \sqrt{1-s^4}} $$
and hence $\gk(s)\to 0$ and $\gtau(s)\to 0$ as $s\to 0$, whereas both $\gk$ and $\gtau$ are summable functions in $L^1(-1,1)$.
Moreover, the Frenet-Serret formulas \eqref{FS} hold true in the open intervals $]-1,0[$ and $]0,1[$.
\par Since $\ttt(0)=2^{-1/2}(1,0,1)$, $\ddot c(0)=0_{\gR^3}$, and $c^{(3)}(0)=2^{-1/2}(0,1,0)$, by using the formulas in \eqref{tbn} we get:
$$ \bbb(0):=\frac{\dot c(0)\tim c^{(3)}(0) }{ \Vert c^{(3)}(0)\Vert }={1\over\sqrt 2}\,( -1,0,1)\,,\quad
\nnn(0):=\bbb(0)\tim\ttt(0)=(0,1,0)
$$
and hence
both the unit normal and binormal are not continuous at $s=0$. However, since $[\nnn(s)]\to[\nnn(0)]$ and $[\bbb(s)]\to[\bbb(0)]$ as $s\to 0$,
they are both continuous as functions with values in $\RP$.
\par We finally compute the total curvature and the total absolute torsion of $c$. With $t=s^2$, we have:
$$ \TC(c)=\int_{-1}^1\gk(s)\,ds=\int_{-1}^1\frac{\sqrt 2 |s|}{ \sqrt{1-s^4}}\,\,ds=\sqrt 2\int_{0}^1\frac{ 1}{\sqrt{1-t^2}}\,dt=\frac\p{\sqrt 2}$$
and similarly
$$ \TAT(c)=\int_{-1}^1|\gtau(s)|\,ds=\int_{-1}^1\frac{\sqrt 2 |s|}{ \sqrt{1-s^4}}\,\,ds=\frac\p{\sqrt 2}\,.$$
In fact, $c$ is regular at $s=0$, so that there is no turning angle at $c(0)$, whereas
$\bbb(0-)=-\bbb(0+)$, so that also the total absolute torsion is zero at $c(0)$.
On the other hand, due to the occurrence of an inflection point at $c(0)$, the complete torsion in the sense of Alexandrov-Reshetnyak \cite{AR} yields a contribution equal to $\p$ at $c(0)$, so that $\CT(c)=\TAT(c)+\p$.
\eex
\br\label{Rsmooth} We finally point out that with the above assumptions, the statements of Theorem~\ref{Tbinsmooth}, Proposition~\ref{Ptansmooth}, and Proposition~\ref{Pnorsmooth} continue to hold.
More precisely, using that the non-negative curvature $\gk(\l)$ and the torsion $\gtau(\l)$ may vanish only at a negligible set of inflection points,
with our previous notation one readily obtains the following relations concerning the trihedral $(\ttt,\bbb,\nnn)$ from Proposition~\ref{PFS}\,:
\ben
\item $\ttt(s_1(k))=\ttt_c(k)\in\SP$ for $k\in[0,C]$, where $s_1:[0,C]\to[0,L]$ is the inverse of the function
\beq\label{ks} k(s):=\int_0^s\gk(\l)\,d\l\,,\qquad s\in[0,L] \,; \eeq
\item $[\bbb(s_2(t))]=\bbb_c(t)\in\RP$ for $t\in[0,T]$, where $s_2:[0,T]\to[0,L]$ is the inverse of the function
\beq\label{taus} t(s):=\int_0^s|\gtau(\l)|\,d\l\,, \qquad s\in[0,L]\,;  \eeq
\item $[\nnn(s_3(\r))]=\nnn_c(\r)\in\RP$ for $\r\in[0,C+T]$, where $s_3:[0,C+T]\to[0,L]$ is the inverse of the function
$$\r(s):=\int_0^s(\gk(\l)+|\gtau(\l)|)\,d\l\,, \qquad s\in[0,L]\,.  $$
\een \er
\bex\label{Eflex2} Going back to Example~\ref{Eflex}, we compute
$$ k(s):=\int_{-1}^s\gk(\l)\,d\l={1\over\sqrt 2}\,\Bigl(\frac\p 2+\frac s{|s|}\,\arcsin (s^2) \Bigr)\,,\quad s\in[-1,1]$$
and hence $s_1(k)=|\cos(\sqrt 2k)|^{1/2}$, where $k\in[0,C]$, with $C=\p/\sqrt 2$, so that
$$ \ttt_c(k):=\ttt(s_1(k))=\frac 1{\sqrt 2}\,\bigl(1,|\cos(\sqrt 2 k)|,\sin(\sqrt 2k) \bigr)\,,\qquad k\in[0,\p/\sqrt 2] $$
with $k(0)=\p/(2\sqrt 2)$ and $\ttt_c(k(0))=2^{-1/2}(1,0,1)$.
Notice moreover that
$$ \dot\ttt_c(k)=\left\{ \ba{ll} \bigl(0,-\sin(\sqrt 2 k), \cos(\sqrt 2k)\bigr) & \If k\in[0,\p/(2\sqrt 2) [ \\
\bigl(0,\sin(\sqrt 2 k), \cos(\sqrt 2k)\bigr) & \If k\in]\p/(2\sqrt 2),\p/\sqrt 2] \ea \right. $$
so that $\dot\ttt_c(k(0)\pm)=(0,\pm 1,0)$.
We also get
$$ \bbb_c(t)=\bigl[2^{-1/2}\,\bigl(-1,|\cos(\sqrt 2 t)|,\sin(\sqrt 2 t) \bigr) \bigr]\,, \qquad k\in[0,T]\,,\quad T=\p/\sqrt 2 $$
where $t(0)=\p/(2\sqrt 2)$ and $\bbb_c(t(0))=\bigl[2^{-1/2}\,(-1,0,1) \bigr]$. Finally,
$$ \dot\bbb_c(t)=\left\{ \ba{ll} \bigl[\bigl(0,-\sin(\sqrt 2 k), \cos(\sqrt 2k)\bigr)\bigr] & \If t\in[0,\p/(2\sqrt 2) [ \\
\bigl[(0,\sin(\sqrt 2 k), \cos(\sqrt 2k)\bigr] & \If t\in]\p/(2\sqrt 2),\p/\sqrt 2] \ea \right. $$
so that $\dot\bbb_c(t(0)+)=\dot\bbb_c(t(0)-)=[(0,1,0)]$, whence $\bbb_c$ has no corner points.
\eex
\section{Torsion force}\label{Sec:force}
The {\em curvature force} was introduced in \cite{CFKSW}, see also \cite{Su_curv}, as the distributional derivative of the tangent indicatrix of curves in $\RN$ with finite total curvature.
It comes into the play when one computes the first variation of the length.
\par More precisely, a rectifiable curve $c$ has finite total curvature if and only if the tantrix $\ttt$ is a function with bounded variation, i.e., the distributional derivative
$D\ttt$ is equal to a finite measure, the curvature force $\calK$.
Also, this property is equivalent to the requirement that the first variation $\d_\x\calL(c)$ of the length has distributional order one.
\par In this section, we shall see that a {\em torsion force} measure can be similarly obtained by means of the tangential variation
of the length $\calL_{\SP}(\ttt_c)$ of the tangent indicatrix $\ttt_c$ that we built up in Proposition~\ref{Ptan} for any curve $c$ with finite total curvature.
\par In fact, the first variation $\d_\x\calL_{\SP}(\ttt_c)$ has distributional order one if and only if
the arc-length derivative $\dot\ttt_c$ of the tantrix $\ttt_c$ is a function of bounded variation, see \eqref{vartan}.
By the way, we recall that this condition is satisfied if in addition the curve $c$ has finite complete torsion, $\CT(c)<\i$, see Remark~\ref{Rtan}.
\par In this case, there exists a finite measure, the torsion force $\calT$, such that
$\lan \calT,\x\ran=\lan D\dot\ttt_c,\x\ran$ for each
smooth tangential vector field $\x$ along $\ttt_c$.
\par Finally, the tangential variation of the length of the weak binormal $\bbb_c$ from Theorem~\ref{Tbin} is briefly discussed.
\adl\par\noindent{\large\sc Curvature force.} Let $c:[0,L]\to\RN$ denote a rectifiable curve parameterized in arc-length $s$.
Suppose that $c_\e$ is a variation of $c$ under which the motion of each point $c(s)$ is smooth in time and with initial velocity $\x(s)$, where $\x:[0,L]\to\RN$ is a Lipschitz continuous function of arc length. The first variation formula gives
$$ \d_\x\calL(c):={d\over d\e}\,\calL(c_\e)_{\vert \e=0}=\int_0^L \ttt(s)\bullet\dot\x(s)\,ds $$
where $\ttt(s)=\dot c(s)$ and $\dot\x(s)$ are defined for a.e. $s$, by Rademacher's theorem.
\par If $c$ is of class $C^2$, integrating by parts one gets
$$ \d_\x\calL(c)=-\int_0^L \dot\ttt(s)\bullet\x(s)\,ds+\bigl(\ttt(L)\bullet\x(L)-\ttt(0)\bullet\x(0)\bigr)$$
where in terms of the (positive) first curvature $\gk$ and first unit normal $\nnn(s)$ one has $\dot\ttt(s)=\gk(s)\,\nnn(s) $, see \eqref{FS} for the case $N=3$.
\par More generally, if $c$ is a curve with finite total curvature, then $\ttt$ is a function of bounded variation,
the right and left limits $\ttt(s\pm):=\ds\lim_{\l\to s^\pm}\ttt(\l)\in\mathbb S^{N-1}$ are well defined for each $s\in]0,L[$,
and the distributional derivative $D\ttt$ is a finite vector-valued measure. Therefore, if in addition $\x(0)=\x(L)=0$ one obtains
$$ \d_\x\calL(c)=\int_0^L \ttt(s)\bullet\dot\x(s)\,ds=-\lan D\ttt,\x\ran\,. $$
\par The measure $\calK:=D\ttt$ is called in \cite{CFKSW} the curvature force, and in the smooth case one has $\calK=\gk\,\nnn\,ds$. If $c$ is a piecewise smooth function,
one has the decomposition $\calK=\calK^a+\calK^s$, where the absolutely continuous component $\calK^a$ is equal to $\gk\,\nnn\,d\calL^1\pri]0,L[$, whereas the singular component $\calK^s$ is given by a sum of Dirac masses concentrated at the corner points of the curve $c$.
\par More precisely, taking for simplicity $N=3$, if $s\in]0,L[$ is such that $\ttt(s-)\neq\ttt(s+)$, then $\calK(\{c(s)\})=(\ttt(s+)-\ttt(s-))\,\d_{c(s)}$.
Therefore, if $\t\in]0,\p]$ is the shortest angle in the Gauss sphere $\SP$ between $\ttt(s\pm)$, so that $d_\SP(\ttt(s+),\ttt(s-))=\t$, one has
$|\calK|(\{c(s)\})=\Vert\ttt(s+)-\ttt(s-)\Vert=2\sin(\t/2)$.
\par As a consequence, compare \cite{Su_curv}, denoting by $\TC^*(c)$ the total variation of the curvature force $\calK$, in general one has $\TC^*(c)\leq\TC(c)$, and the strict inequality holds true as soon has the curve $c$ has an interior corner point. More precisely, by the previous computation one has
$$ \TC^*(c)=\Var_{\gR^3}(\ttt)\,,\qquad \TC(c)=\Var_{\SP}(\ttt)\,. $$
{\large\sc First variation of total curvature.} In Proposition~\ref{Ptan},
for any curve $c$ in $\gR^3$ with finite total curvature $C:=\TC(c)$ we have constructed a rectifiable curve $\ttt_c:[0,C]\to\SP$ parameterized in arc-length, that is strictly related with the complete tangent indicatrix in the sense of \cite{AR}.
We wish to compute the tangential variation of the length $\calL_{\SP}(\ttt_c)$ of $\ttt_c$, by considering in particular the smooth case.
\par For this purpose, we assume that $\ttt_{c,\,\e}$ is a variation of $\ttt_c$ under which the motion of each point $\ttt_c(k)$ is smooth in time and with initial velocity $\x(s)$, where this time $\x:[0,C]\to\gR^3$ is a Lipschitz continuous function of arc length $k$, with $\x(0)=\x(C)=0$.
Since we deal with tangential variations, {\em we assume in addition that $\x(k)\in T_{\ttt_c(k)}\SP$ for each $k$}.
The first variation formula gives:
$$ \d_\x\calL_{\SP}(\ttt_c):={d\over d\e}\,\calL_\SP(\ttt_{c,\,\e})_{\vert \e=0}=\int_0^C \dot\ttt_c(k)\bullet\dot\x(k)\,dk $$
where $\dot \ttt_c(k)$ and $\dot\x(k)$ are defined for a.e. $k$.
Therefore, by the definition of distributional derivative, in general we obtain:
\beq\label{vartan} \d_\x\calL_{\SP}(\ttt_c)=\int_0^C \dot\ttt_c(k)\bullet\dot\x(k)\,dk=:-\lan D\dot\ttt_c,\x\ran\,. \eeq
%
%Moreover, from the proof of Proposition~\ref{Ptan} it follows that
%
\par Assume now that $c$ is of class $C^3$ and $\ddot c(s)\neq 0_{\gR^3}$ for each $s\in]0,L[$. In point~i) of Remark~\ref{Rsmooth}, we have seen that
$\ttt_c(k)=\ttt(s_1(k))$ for each $k\in[0,C]$, where $\ttt(s)=\dot c(s)$ and $s_1:[0,C]\to[0,L]$ is the inverse of the function $k(s)$ in \eqref{ks},
so that
$$\dot s_1(k)=\gk(s_1(k))^{-1}\,,\qquad \dot\ttt_c(k)=\ttt'(s_1(k))\, \dot s_1(k)=\nnn(s_1(k)) $$
for each $k\in[0,C]$, by the first Frenet-Serret formula in \eqref{FS}. Therefore,
by the second formula in \eqref{FS} we compute for each $k$
$$ \ddot \ttt_c(k)=\nnn'(s_1)\, \dot s_1(k)=  -\ttt(s_1)+{\gtau(s_1)\over \gk(s_1)} \,\bbb(s_1)\,,\qquad s_1=s_1(k)\,.$$
\par Now, the tangential component to $\SP$ of the second derivative $\ddot \ttt_c(k)$, i.e.,
the {\em geodesic curvature} of $\ttt_c$ at the point $\ttt_c(k)$, agrees with the quotient between the torsion and the scalar curvature of the initial curve $c$ at the point $c(s_1)$, where $s_1=s_1(k)$.
\par In fact, the Darboux frame along $\ttt_c$ is the triad $(\gT,\gN,\gU)$,
where $\gT(k):=\dot\ttt_c(k)$, $\gN(k):=\nu(\ttt_c(k))$, $\nu(p)$ being the unit normal to the tangent 2-space $T_p\SP$,
and $\gU(k):=\gN(k)\tim\gT(k)$ is the unit conormal.
The curvature vector $\gK(k):=\dot\gT(k)=\ddot\ttt_c(k)$ is orthogonal to
$\gT(k)$, and thus decomposes as
$$\gK(k)=\kk_g(k)\,\gU(k)+\kk_n(k)\,\gN(k) $$
where $\kk_g:=\gK\bullet\gU$ and $\kk_n:=\gK\bullet\gN$ denote the geodesic and normal curvature of $\ttt_c$, respectively.
By changing variable, we get
$$ \gT(k)=\nnn(s_1)\,,\quad \gN(k)=\ttt(s_1)\,,\quad \gU(k)=\bbb(s_1) $$
and hence we obtain for each $k\in[0,C]$
$$ \kk_g(k)= {\gtau(s_1)\over \gk(s_1)}\,,\qquad \kk_n(k)\equiv -1\,,\qquad s_1=s_1(k)\,. $$
\par As a consequence, integrating by parts in \eqref{vartan} we get
$$ \lan D\dot\ttt_c,\x\ran=\int_0^C \kk_g(k)\,\bbb(s_1(k))\bullet\x(k)\,dk= \int_0^C {\gtau(s_1)\over \gk(s_1)}\,\bbb(s_1)\bullet\x(k)\,dk $$
where, we recall, $\x(k)\in T_{\ttt_c(k)}\SP$ for each $k$. Therefore, by changing variable $s=s_1(k)$, since $ds=\gk(s_1(k))^{-1}\,dk$ we recover the expected formula:
$$ %\d_\x\calL_{\SP}(\ttt_c)
\lan D\dot\ttt_c,\x\ran= \int_0^L \gtau(s)\,\bbb(s)\bullet\x(k(s))\,ds\,. $$
{\large\sc Torsion force.} Denoting by ${\calT}$ the ``tangential" component of the distributional derivative of $\dot\ttt_c$, so that $\lan \calT,\x\ran=\lan D\dot\ttt_c,\x\ran$ for each
smooth tangential vector field along $\ttt_c$, we have just seen that if $c$ is smooth, then
%
%$$ \calT =\gtau\circ s_1\cdot\bbb\circ s_1\,d\calL^1\pri[0,C]\,,\qquad s_1=s_1(k)\,. $$
%
\beq\label{calT}
 \calT = k_{\#}\bigl(\gtau\,\bbb\,d\calL^1\pri]0,L[\bigr) \eeq
i.e., $\cal T$ is the push forward of the measure $\gtau\,\bbb\,d\calL^1\pri]0,L[$ by the function $k(s)$ defined in \eqref{ks}, and its total mass is equal to $\int_c|\gtau|\,ds$.
For that reason, the measure ${\calT}$ may be called the {\em torsion force}.
\par More generally, it turns out that the torsion force $\calT$ is a finite measure provided that the derivative $\dot\ttt_c$ is a function of bounded variation.
To this purpose, see Remark~\ref{Rtan}, we recall that this sufficient condition is satisfied if in addition the curve $c$ has finite complete torsion, $\CT(c)<\i$.
\par If $c$ is piecewise smooth, we obtain again the decomposition $\calT=\calT^a+\calT^s$. By Proposition~\ref{Ptansmooth}, it turns out that the absolutely continuous component
$ \calT^a $ takes the same form as in the right-hand side of the formula \eqref{calT}, where this time $k(s):=\TC(c_{\vert[0,s]})$.
Moreover, using that  $\ttt(s)=\ttt_c(k(s))$, if $c$ is smooth at $s$ we have $\ttt'(s)=\dot\ttt_c(k(s))\cdot k'(s)$, with $k'(s)=\gk(s)$, hence by the first formula in \eqref{FS} we get $\dot\ttt_c(k(s))=\nnn(s)$.
\par If $c$ has a point of return at $c(s)$, we have $\ttt(s-)=-\ttt(s+)$.
In this case, see Remark~\ref{Runique}, the curve $\ttt_c$ (and hence the torsion force $\calT$) depends on the choice of the geodesic arc connecting the antipodal points $\ttt(s\pm)$.
However, the total mass of $\calT$ is finite and it does not depend on the choice of the geodesics.
\par If $c$ has no points of return, the torsion force $\calT$ only depends on $c$. In fact, the singular component $\calT^s$ is a sum of Dirac masses concentrated at the corner points $x=\ttt_c(k)$ of the curve $\ttt_c$, with weight $\dot\ttt_c(k+)-\dot\ttt_c(k-)$.
If $\t$ is the turning angle of $\ttt_c$ at $x$, then $\Vert\dot\ttt_c(k+)-\dot\ttt_c(k-)\Vert=2\sin(\t/2)$.
\par In Example~\ref{Eflex2}, at $x=\ttt_c(k(0))=2^{-1/2}(1,0,1)$ we have $\dot\ttt_c(k(0)\pm)=(0,\pm 1,0)$, so that $\t=\p$ and
$\Vert\dot\ttt_c(k(0)+)-\dot\ttt_c(k(0)-)\Vert=2$.
\adl\par\noindent
{\large\sc First variation of total torsion.}
In Theorem~\ref{Tbin}, we defined the weak binormal $\bbb_c:[0,T]\to\RP$, that satisfies $|\dot\bbb_c|=1$ a.e. and $\calL_{\RP}(\bbb_c)=\TAT(c)$, and it turns out that the derivative $\dot\bbb_c$ is a function of bounded variation.
\par Moreover, in point ii) of Remark~\ref{Rsmooth}, we have seen that if $c$ is smooth as above, then
$\bbb_c(t)=[\bbb(s_2(t))]$ for each $t\in[0,T]$, where $s_2:[0,T]\to[0,L]$, with $T=\TAT(c)$, is the inverse of the function $t(s)$ in \eqref{taus}.
We have
$$\dot s_2(t)=|\gtau(s_2(t))|^{-1}\,,\qquad \dot\bbb_c(t)=[\bbb'(s_2(t))]\,\dot s_2(t)=-\sgn\bigl( \gtau(s_2(t))\bigr)\,[\nnn(s_2(t))] $$
for each $t$, by the third Frenet-Serret formula in \eqref{FS}. Therefore,
by the second formula in \eqref{FS} we get
$$ {d^2\over dt^2}\bbb(s_2)= -\sgn\bigl( \gtau(s_2)\bigr)\,\nnn'(s_2)\dot s_2(t)= {\gk(s_2)\over \gtau(s_2)}\,\ttt(s_2)-\bbb(s_2) \,,\qquad s_2=s_2(t)\,.$$
%
%$$ \ddot \bbb_c(t)= -\sgn\bigl( \gtau(s_2)\bigr)\,[\nnn'(s_2)]\dot s_2(t)= {\gk(s_2)\over \gtau(s_2)}\,[\ttt(s_2)]-[\bbb(s_2)] \,,\qquad s_2=s_2(t)\,.$$
%
\par Arguing as above, the tangential variation of the length $\calL_{\RP}(\bbb_c)$ yields to the ``tangential" component of the distributional derivative
$D\dot\bbb_c$. In the smooth case, its lifting gives the measure
$t_{\#}\bigl(\sgn(\gtau)\gk\,\nnn\,d\calL^1\pri]0,L[\bigr)$, with total mass $\int_c|\gk|\,ds$.
If $c$ is piecewise smooth, the singular component $\calT^s$ is a sum of Dirac masses concentrated at the corner points $x=\bbb_c(t)$ of the curve $\bbb_c$ in $\RP$, i.e., at the points where $\dot\bbb_c(t+)\neq\dot\bbb_c(t-)$ in $\RP$.
Notice however that in Example~\ref{Eflex2}, at $x=\bbb_c(t(0))=\bigl[2^{-1/2}\,(-1,0,1) \bigr]$ we have
$\dot\bbb_c(t(0)+)=\dot\bbb_c(t(0)-)=[(0,1,0)]$, whence $\bbb_c$ has no corner points and the measure $D\dot\bbb_c$ has no singular part.
%
%For that reason, up to a lifting one essentially recovers the curvature force $\calK$ from \cite{CFKSW}.
%
%
%PRIN 2010-2011 ``Variet\`a reali e complesse: geometria, to\-po\-lo\-gia e analisi ar\-mo\-ni\-ca''
%FIRB 2008 ``Geometria Differenziale Complessa e Dinamica Olomorfa'';
%
\adl\par\noindent{\bf Acknowledgements.}
D.M. wishes to thank his seventeen-year-old son, Gio\-van\-ni Mucci, for his useful help in pointing out to us some interesting phenomena concerning spherical geometry that he learned by himself.
We also thank the referee for his/her suggestions that allowed us to improve the first version of the paper and to obtain the results contained in the last section.
The research of D.M. was partially supported by the GNAMPA of INDAM.
The research of A.S. was partially supported by the GNSAGA of INDAM.
\end{document}